\documentclass[reqno,11pt]{amsart}
\usepackage{mathrsfs,graphicx,amsthm,amsfonts,amsxtra,xspace,array,cite,epic,float,ifpdf,bm}
\usepackage[left=2.5cm]{geometry}
\usepackage{url}
\usepackage{enumerate}
\usepackage{asymptote}
\usepackage{amsmath,amssymb,amscd}
\usepackage{ifthen}     %%%%%   to enable swithcing between choices %%%

\newcommand\pubpri[2]{%
\ifthenelse{\equal{\version}{public}}%
{{#1}}%
{\marginpar{\scshape\small Pubpri Alert}{#2}}}
\newcommand\pubprinoalert[2]{%
\ifthenelse{\equal{\version}{public}}%
{{#1}}%
{#2}}
\newcommand\ignore[1]{}
\usepackage{color}
\providecommand\wantcolor{yes}   %  
\definecolor{backgroundyellow}{cmyk}{.2,.1,.8,.2}
\definecolor{backgroundblue}{rgb}{0,0,1}
\definecolor{backgroundred}{rgb}{1,0,0}
\definecolor{backgroundmagenta}{cmyk}{0,1,0,0}
\ifthenelse{\equal{\wantcolor}{no}}{\renewcommand\color[1]{}}{}

\newcommand\mysubsubsection[1]{%
		\subsubsection{\sffamily\upshape\mdseries #1}}

\newcommand\mysss{\mysubsubsection}

\ifpdf
  \usepackage[
    pdftex,
    colorlinks,%
    linkcolor=blue,citecolor=red,urlcolor=blue,
    hyperindex,%
    plainpages=false,%
    bookmarksopen,%
    bookmarksnumbered%
  ]{hyperref} 
\else
  \usepackage{hyperref}
\fi
\newtheorem{annotation}{Annotation}[section]
\newtheorem{theorem}[annotation]{%\color{blue}
		Theorem}
\newtheorem{lemma}[annotation]{%\color{blue}
		Lemma}
\newtheorem{definition}[annotation]{%\color{magenta}
		Definition}
\newtheorem{corollary}[annotation]{%\color{blue}
		Corollary}
\newtheorem{proposition}[annotation]{%\color{cyan}
		Proposition}
\newtheorem{example}[annotation]{%\color{yellow}
		Example}
\newtheorem{conjecture}[annotation]{Conjecture}		
\newcommand\bexample{\begin{example}\begin{rm}}
\newcommand\eexample{\end{rm}\hfill$\Box$\end{example}}
\newtheorem{examplenobox}[annotation]{%\color{yellow}
		Example}
\newcommand\bexamplenobox{\begin{examplenobox}\begin{rm}}
\newcommand\eexamplenobox{\end{rm}\end{examplenobox}}
\newtheorem{exercise}[annotation]{%\color{yellow}
		Exercise}
\newcommand\bexercise{\begin{exercise}\begin{rm}}
\newcommand\eexercise{\end{rm}\end{exercise}}
\newtheorem{notation}[annotation]{%\color{yellow}
		Notation}
\newcommand\bnotation{\begin{notation}\begin{rm}}
\newcommand\enotation{\end{rm}\end{notation}}

\newtheorem{remark}[annotation]{%\color{yellow}
		Remark}
\newcommand\bremark{\begin{remark}%\begin{mdseries}\begin{sffamily}
\begin{upshape}}
\newcommand\eremark{\end{upshape}%\end{sffamily}\end{mdseries}
\end{remark}}
\newcommand\bdefn{\begin{definition}%\begin{mdseries}\begin{sffamily}
\begin{upshape}}
\newcommand\edefn{\end{upshape}%\end{sffamily}\end{mdseries}
\end{definition}}
\newtheorem{caveat}[annotation]{%\color{yellow}
		Caveat}
\newcommand\bcaveat{\begin{caveat}%\begin{mdseries}\begin{sffamily}
\begin{upshape}}
\newcommand\ecaveat{\end{upshape}%\end{sffamily}\end{mdseries}
\end{caveat}}
\newenvironment{caveatstar}{%\color{yellow}
\par\noindent{\scshape\bfseries
  Caveat: }\begin{rm}}{\end{rm}\newline} 
\newcommand\bcaveatstar{\begin{caveatstar}}%\begin{mdseries}\begin{sffamily}
\newcommand\ecaveatstar{\end{caveatstar}}
\newenvironment{myproof}{%
\par\noindent{\scshape
  Proof: }\begin{rm}}{\hfill$\Box$\end{rm}\newline} 
\newcommand\bmyproof{\begin{myproof}}
\newcommand\emyproof{\end{myproof}}
\newenvironment{myproofnobox}{%
\par\noindent{\scshape Proof:}\begin{rm}}{\end{rm}\hfill\newline}
\newcommand\bmyproofnobox{\begin{myproofnobox}}
\newcommand\emyproofnobox{\end{myproofnobox}}
\newenvironment{solution}{%
\par\noindent{\scshape Solution: }\begin{rm}}{\hfill$\Box$\end{rm}\newline}
\newenvironment{solutionnobox}{%
\par\noindent{\scshape Solution: }\begin{rm}}{\end{rm}}
\newcommand\bsolution{\begin{solution}\begin{rm}}
\newcommand\esolution{\end{rm}\end{solution}}
\newcommand\bsolutionnobox{\begin{solutionnobox}\begin{rm}}
\newcommand\esolutionnobox{\end{rm}\end{solutionnobox}}
\newcommand\bthm{\begin{theorem}}
\newcommand\ethm{\end{theorem}}
\newcommand\bcor{\begin{corollary}}
\newcommand\ecor{\end{corollary}}
\newcommand\blemma{\begin{lemma}}
\newcommand\elemma{\end{lemma}}
\newcommand\bprop{\begin{proposition}}
\newcommand\eprop{\end{proposition}}
\newcommand\beqn{\begin{equation}}
\newcommand\eeqn{\end{equation}}
\newcommand\beqnstar{\begin{equation*}}
\newcommand\eeqnstar{\end{equation*}}
\newcommand\mtitle[1]%
{\marginpar{\sffamily\raggedright\upshape\small #1}}

\providecommand\finalized{no}
\ifthenelse{\equal{\finalized}{yes}}{}{\marginparwidth=2cm}
\ifthenelse{\equal{\finalized}{yes}}%
{}%
{}
\ifthenelse{\equal{\finalized}{yes}}%
{\newcommand\checked[1]{}}%
{\newcommand\checked[1]{\marginpar{[{\ttfamily\upshape\tiny CHECKED: #1}]}}}
\ifthenelse{\equal{\finalized}{yes}}%
{\newcommand\spellchecked[1]{}}%
{\newcommand\spellchecked[1]{\marginpar{[{\ttfamily\upshape\tiny SPELLCHECKED: #1}]}}}
\providecommand\version{public}   %  setting default  
\ifthenelse{\equal{\version}{public}}%
{\newcommand\mcomment[1]{}}%
{\newcommand\mcomment[1]{\marginpar{{\raggedright\sffamily\upshape\small
\begin{spacing}{0.75} #1\end{spacing}}}}}
\ifthenelse{\equal{\version}{public}}%
{\newcommand\fcomment[1]{}}%
{\newcommand\fcomment[1]{\footnote{#1}}}
\ifthenelse{\equal{\version}{public}}%
{\newcommand\comment[1]{}}%
{\newcommand\comment[1]{{\small #1}}}
\renewcommand\Bbb{\mathbb}

\newcommand{\be}{\begin{enumerate}}
\newcommand{\ee}{\end{enumerate}}
\newcommand{\beq}{\begin{equation}}
\newcommand{\eeq}{\end{equation}}
\newcommand{\beqs}{\begin{equation*}}                     
\newcommand{\eeqs}{\end{equation*}}
\newcommand{\complex}{\mathbb{C}}

\renewcommand\omega{\varpi}

\newcommand{\csa}{\mathfrak{h}}

\newcommand{\csaf}{\widehat{\mathfrak{h}}}
\newcommand\csahat{\widehat{\mathfrak{h}}}

\newcommand{\afw}{\widehat{W}} % the affine Weyl group
\newcommand{\digautos}{\Sigma} %dynkin diag autos

\newcommand\lieg{\mathfrak{g}}

\newcommand{\nminus}{\mathfrak{n}^-}
\newcommand{\nplus}{\mathfrak{n}^+}

\newcommand{\kay}{K}
\newcommand{\hi}{\alpha^\vee_i}
\newcommand{\iprime}{i^\prime}

\newcommand{\bform}[2]{\left(#1 | #2\right)}

\newcommand{\lambar}{\lambda}

\newcommand{\lamprbar}{{\lambda^\prime}}
\title[Generalized Demazure modules and fusion products]{Generalized Demazure modules and fusion products}
\author{B.~Ravinder}
\address{Tata Institute of Fundamental Research\\
Homi Bhabha Road, Colaba\\
Mumbai 400005, India}
\email{ravinder@math.tifr.res.in}

\begingroup\makeatletter\ifx\SetFigFont\undefined%
\gdef\SetFigFont#1#2#3#4#5{
  \reset@font\fontsize{#1}{#2pt}
  \fontfamily{#3}\fontseries{#4}\fontshape{#5}
  \selectfont}
\fi\endgroup
\begin{document}
\allowdisplaybreaks
\numberwithin{equation}{section}
\begin{abstract}
Let $\mathfrak{g}$ be a finite-dimensional complex simple Lie algebra with highest root $\theta$ and let $\mathfrak{g}[t]$ be the corresponding current algebra. 
In this paper, we  consider the $\mathfrak{g}[t]$-stable Demazure modules associated to integrable highest weight representations of the affine Lie algebra $\widehat{\mathfrak{g}}$. 
We prove that the fusion product of Demazure modules of a given level with a single Demazure module of a different level and with highest weight a multiple of $\theta$ is  
a generalized Demazure module, and also give defining relations. This also shows that the fusion product of such Demazure modules is independent of the chosen parameters. 
As a consequence we obtain generators and relations for certain types of generalized Demazure modules. We also establish a connection with the modules defined 
by Chari and Venkatesh.
\end{abstract}

\maketitle
\section{Introduction}
Let $\lieg$ be a finite-dimensional complex simple Lie algebra and $\widehat{\lieg}$ the corresponding affine Lie algebra. In this paper, we are interested in Demazure modules associated to integrable highest weight representations of $\widehat{\lieg}$.
These modules, which are actually modules  for a Borel subalgebra of $\widehat{\lieg}$, are indexed by a dominant integral affine weight  and an element of the affine Weyl group. We are mainly interested in the Demazure modules which are preserved by a maximal parabolic subalgebra containing the Borel. 
The maximal parabolic subalgebra of our interest is the current algebra $\lieg[t]$, which is the algebra of polynomial maps $\complex\to\lieg$ with the obvious 
point-wise bracket. Equivalently, it is the complex vector space $\lieg\otimes\complex[t]$ with Lie bracket being the $\complex[t]$-bilinear extension of the Lie bracket 
on $\lieg$. The degree grading on $\complex[t]$ gives a natural $\mathbb{Z}$-grading on $\lieg[t]$ and makes it a graded Lie algebra. The $\lieg[t]$-stable Demazure modules are known to be indexed  by pairs $(\ell, \lambda)$, 
where $\ell$ is the level of the integrable representation of $\widehat{\lieg}$ and $\lambda$ is a dominant integral weight of $\lieg$. 
We denote the corresponding module by $D(\ell,\,\lambda)$. These are in fact finite-dimensional graded $\lieg[t]$-modules.

A powerful tool to study the category of finite-dimensional graded $\lieg[t]$-modules is the fusion product, which was introduced by Feigin and Loktev in \cite{FL}. Although the fusion product is by definition dependent on a choice of parameters, it is widely expected that it will turn out to be independent of the choices, and in several cases this has been proved (see \cite{CSVW, CV, FoL, KL, Naoi1, Naoi2, R1}).
It is proved in \cite{CSVW} that the fusion product of Demazure modules of a given level is again a Demazure module of the same level. In \cite{Naoi2}, for $\lieg$ simply laced, it is proved that the fusion product of Demazure modules of different level with highest weight a multiple of a fundamental weight is a generalized Demazure module, and used this to solve the $X=M$ conjecture.
The generalized Demazure modules are indexed by $p$ dominant integral affine weights and $p$ affine Weyl group elements, where $p\geq 1$.  Their defining relations are not known except when $p=1$, where they are actually the Demazure modules. But in special cases the character of generalized Demazure modules  is known in terms of the Demazure operators \cite{LLM, Naoi3}. 

In this paper, we investigate further the fusion product of different level Demazure modules. Let $\theta$ be the  highest root  of $\lieg$. We consider the fusion product of Demazure modules of a given level with a single Demazure module of a different level and with highest weight a multiple of $\theta$. We prove that this fusion product 
as a $\lieg[t]$-module is isomorphic to a generalized Demazure module, and give the defining relations. 
More precisely, given positive integers $k,\ell,m$ such that $\ell\geq m\geq k$, and a sequence of dominant integral weights $\lambda_1,\ldots,\lambda_p$ of $\lieg$,  we prove that the fusion product 
\beq\label{fus}D(\ell,\,\ell\lambda_1)*\cdots * D(\ell,\,\ell\lambda_p)*D(m,\,k\theta),\eeq of Demazure modules 
is a generalized Demazure module, and also give the defining relations (see Theorem \ref{MTtwo}). This also proves that the fusion product \eqref{fus} is independent of the chosen parameters. 

Our main results (Theorems \ref{MTone} and \ref{MTtwo}) enable us to obtain short exact sequences of fusion products and  generalized Demazure modules (see Corollary \ref{exactfg}), and a surjective morphism between two fusion products (see Corollary \ref{Schur}). As a consequence of 
Corollary \ref{Schur}, we get the following result which may be viewed as a generalization of the Schur positivity \cite{CFS} (see Corollary \ref{Schurpm}):
Given two partitions $(\ell_1\geq \ell_2\geq\cdots\geq \ell_{p}\geq0)$ and $(m_1\geq m_2\geq\cdots\geq m_{p}\geq0)$ of a positive integer, 
%such that $\bm{m}\preceq{\bm\ell}$, where $\preceq$
%is the reverse dominance order on partitions of $N$, 
there exists a surjective morphism of $\lieg$-modules
$$D(\ell_1,\,\ell_1\theta)\otimes D(\ell_2,\,\ell_2\theta)\otimes\cdots\otimes D(\ell_{p},\,\ell_{p}\theta)\twoheadrightarrow D(m_1,\,m_1\theta)\otimes D(m_2,\,m_2\theta)\otimes\cdots\otimes D(m_{p},\,m_{p}\theta),$$
if $\ell_i+\cdots+\ell_p\geq m_i+\cdots+m_p$ holds for each $1\leq i\leq p$.

In \cite{CV}, Chari and Venkatesh have introduced a large family of indecomposable graded $\lieg[t]$-modules which includes the Demazure modules. 
In \S\ref{connectiontoCV}, we prove that certain types of generalized Demazure modules also belong to their family (see Corollary \ref{cor:all}).

\subsection*{Acknowledgements}
The author thanks K.~N.~Raghavan and S.~Viswanath for many helpful discussions and encouragement. The author acknowledges support from TIFR, Mumbai, under the Visiting Fellowship scheme. Most of this work was done when the author
was a Ph.D. student at IMSc, Chennai, and he acknowledges support from CSIR under the SPM Fellowship scheme.
\section{Preliminaries}
Throughout the paper, $\mathbb{Z}$ denotes the set of integers, $\mathbb{N}$  the set of positive integers,  $\mathbb{Z}_{\geq 0}$ the set of non-negative integers,   $\complex$  the field of complex numbers, $\complex[t]$ the polynomial ring in an indeterminate $t$ and  $\complex[t,t^{-1}]$ the ring of Laurent polynomials. 
\subsection{}
Given a complex Lie algebra $\mathfrak{a}$, let $\mathbf{U}(\mathfrak{a})$ be its universal enveloping algebra. The current algebra $\mathfrak{a}[t]$  associated to $\mathfrak{a}$ is defined as $\mathfrak{a} \otimes \complex[t]$, with the Lie bracket
\[ [x \otimes t^r,  y \otimes t^s]=[x, y]\otimes t^{r+s} \quad \forall \,\, x, y \in  \mathfrak{a},\,\, r, s \in \,\mathbb{Z}_{\geq 0}.\]
The degree grading on $\complex[t]$ gives a natural $\mathbb{Z}_{\geq 0}$-grading on $\mathbf{U}(\mathfrak{a}[t])$: the element $(a_1\otimes t^{s_1})\cdots (a_k\otimes t^{s_k})$, for $a_i\in\mathfrak{a}, s_i\in \mathbb{Z}_{\geq 0},$ has grade $s_1+\cdots+s_k.$
A graded $\mathfrak{a}[t]$-module is a $\mathbb{Z}$-graded vector space
$V=\bigoplus_{r\in\mathbb{Z}} V[r]$ such that
\[ (\mathfrak{a}\otimes t^s) V[r]\subset V[r+s], \quad \forall\,\, r\in\mathbb{Z},\, s\in\mathbb{Z}_{\geq0}.\]
%Given a graded $\mathfrak{a}[t]$-module $V$ and a subset $S$ of $V$, $\langle S \rangle$ denotes the submodule of $V$ generated by $S$. 
Let $\textup{ev}_0:\mathfrak{a}[t]\rightarrow \mathfrak{a}$ be the morphism of Lie algebras given by setting $t=0$. The pull back of any $\mathfrak{a}$-module $V$ by $\textup{ev}_0$ defines a graded $\mathfrak{a}[t]$-module structure on $V,$ and we denote this module by $\textup{ev}_0\,V.$
We define the morphism of graded $\mathfrak{a}[t]$-modules as a degree zero morphism of $\mathfrak{a}[t]$-modules. For $r\in\mathbb{Z}$ and a graded $\mathfrak{a}[t]$-module $V,$ we let $\tau_r V$ be the $r$-th graded shift of $V.$ 
\subsection{}
Let $\lieg$ be a finite-dimensional simple Lie algebra over $\complex$ of rank $n$. 
Fix a Cartan subalgebra $\csa$ of $\lieg$ and a Borel subalgebra $\mathfrak{b}$  of $\lieg$ containing $\csa$. 
Let $\lieg=\nminus\oplus\csa\oplus\nplus$ be the triangular decomposition of $\lieg$ with $\mathfrak{b}=\csa\oplus\nplus$.
Let $R \, (\text{resp}. \, R^+)$ be the set of roots (resp. positive roots) of $\lieg$. Let $\theta \in R^+$ be the highest root of $\lieg$. 
Let $(.|.)$ be a non-degenerate, symmetric, invariant  bilinear form  on $\csa^*$ normalized so that the square length of a long root is two. 
It is easy to see from the abstract theory of root systems that $(\theta|\alpha)\in\{0,1\}$, for all $\alpha\in R^+\setminus\{\theta\}$.
For $\alpha\in R,$ let $\alpha^\vee\in\csa$ denote the corresponding co-root, $\lieg_{\alpha}$ the corresponding root space of $\lieg$, and we fix non-zero elements $x^{\pm}_{\alpha}\in \lieg_{\pm\alpha}$ such that $[x^{+}_{\alpha}, x^{-}_{\alpha}]=\alpha^{\vee}$. 
For a root $\alpha$, we set $d_{\alpha}=\frac{2}{(\alpha|\alpha)}$, then we have 
$$ d_\alpha=\begin{cases}
             1 & \textup{if} \,\,\alpha \,\,\textup{is long},\\
             2 & \textup{if} \,\,\alpha \,\, \textup{is short and}\,\,\lieg \,\,\textup{is of type}\,\, B_n, C_n \,\,\textup{or}\,\, F_4,\\
             3 & \textup{if} \,\,\alpha \,\, \textup{is short and}\,\,\lieg \,\,\textup{is of type}\,\, G_2.
            \end{cases}
$$

Set $I=\{1,2,\ldots, n\}$. Let $\alpha_i$ and $\omega_i, i\in I$, be simple roots and fundamental weights respectively. For $\alpha=\sum_{i}n_i\alpha_i\in R$,
we define the {\em height} of $\alpha$ by $\textup{ht}\,\alpha= \sum_{i}n_i$.
For $i\in I$, we write $x^{\pm}_i$ for $x^\pm_{\alpha_i}$ and $d_i$ for $d_{\alpha_i}$.
The weight lattice $P$ (resp. $P^+$) is the $\mathbb{Z}$-span (resp. $\mathbb{Z}_{\geq 0}$-span) of  $\{\omega_i:i\in I\}$. The root lattice $Q$ 
(resp. $Q^+$) is the $\mathbb{Z}$-span (resp. $\mathbb{Z}_{\geq 0}$-span) of $\{\alpha_i:i\in I\}$. The co-weight lattice $L=\sum_{i\in I}\mathbb{Z}d_{i}\omega_i$ is a 
sub lattice of $P$ and the co-root lattice $M=\sum_{i\in I}\mathbb{Z}d_{i}\alpha_i $ is a sub lattice of $Q$. The subsets $L^+$and $M^+$ are defined in the obvious way. 
It is easy to see that for $\lambda\in L^+$ and $\alpha\in R^+$, we have $(\lambda|\alpha)\in\mathbb{Z}_{\geq 0}$.
Let $W$ be the Weyl group of $\lieg$. For $\alpha\in R^+$, we denote by $r_\alpha\in W$ the reflection associated with $\alpha$, and set $r_i=r_{\alpha_i} (i\in I)$. We denote by $w_0$  the longest element in $W$.

For $\lambda\in P^+,$ let $V(\lambda)$ be the  corresponding finite-dimensional irreducible $\lieg$-module generated by an element $v_{\lambda}$ with  the following defining relations:
\[ x^{+}_{i}\,v_\lambda=0,\quad \hi\,v_\lambda = \langle \lambda, \, \hi \rangle\, v_\lambda, \quad (x^{-}_{i})^{\langle \lambda, \, \hi \rangle +1} \,v_\lambda=0, \quad \forall \,\, i\in I.\]
\subsection{}
Let $\widehat{\lieg}$ be the affine Lie algebra defined by
\[\widehat{\lieg}=\lieg\otimes\complex[t,t^{-1}]\oplus\complex c\oplus \complex d,\]
where $c$ is central and the other Lie brackets are given by
\[[x\otimes t^r,y\otimes t^s]=[x,y]\otimes t^{r+s}+r\delta_{r,-s}(x|y)c,\]
\[[d,x\otimes t^r]=r(x\otimes t^r),\]
for all $x,y\in\lieg$ and integers $r,s$. 
The Lie subalgebras $\widehat{\csa}$ and $\widehat{\mathfrak{b}}$ of $\widehat{\lieg}$ are defined as follows:
%\[\widehat{\mathfrak{n}}^-=\mathfrak{n}^-\otimes\complex[t^{-1}]\oplus
 % (\mathfrak{n}^+\oplus\csa)\otimes t^{-1}\complex[t^{-1}],\qquad
\begin{gather*}
\widehat{\csa}=\csa\oplus\complex c \oplus \complex d,
\qquad\widehat{\mathfrak{b}}=\lieg\otimes t\complex[t]\oplus\mathfrak{b}\oplus\complex c \oplus \complex d.
\end{gather*}
%The Cartan subalgebra $\widehat{\csa}$ of $\widehat{\lieg}$ is
%\[\widehat{\csa}=\csa\oplus\complex c \oplus \complex d,\]
We regard ${\csa}^{*}$ as a subspace of ${\widehat{\csa}}^{*}$ by setting $\langle\lambda, c\rangle=\langle\lambda, d\rangle=0$ for $\lambda\in{\csa}^{*}$. 
For $\xi\in{\widehat{\csa}}^{*}$, let $\xi|_{\csa}$ be the element of  ${\csa}^{*}$ obtained by restricting $\xi$ to $\csa$.
Let $\delta, \Lambda_0\in{\csahat}^{*}$ be given by 
\[\langle\delta, \csa+\complex c\rangle=0, \, \, \langle\delta, d\rangle=1, \quad \langle\Lambda_0, \csa+\complex d\rangle=0,\,\,\langle\Lambda_0, c\rangle=1.\]
Extend the non-degenerate form on ${\csa}^{*}$ to a non-degenerate symmetric  
bilinear form on ${\csahat}^{*}$ by setting, 
\[(\csa^{*}|\complex \delta+\complex \Lambda_0)=(\delta|\delta)=(\Lambda_0|\Lambda_0)=0\,\, \text{and} \, \, (\delta|\Lambda_0)=1.\]  

Set $\widehat{I}=I\cup\{0\}$.
For $i\in I$, denote $\Lambda_i=\omega_i+\langle \omega_i,\,\theta^\vee \rangle \Lambda_0 \in \widehat{\csa}^*$.
Let ${\widehat{P}}^+=\sum_{i\in\widehat{I}} \mathbb{Z}_{\geq 0}(\Lambda_i+\mathbb{Z}\delta)$ be the set of dominant integral affine weights, and $\widehat{P}$ is defined similarly. 
%The affine root lattice $\widehat{Q}$ is the $\mathbb{Z}$-span of the simple roots $\alpha_i, i\in\widehat{I}$, of $\widehat{\lieg}$, and $\widehat{Q}^+$ is defined in the obvious way. Let $\widehat{R}_{re}=\{\alpha+r\delta : \alpha \in R,\,r\in\mathbb{Z}\}$ be the set of real roots, $\widehat{R}_{im}=\{r\delta : r\in\mathbb{Z}\setminus \{0\}\}$ the set of imaginary roots and $ \widehat{R}=\widehat{R}_{re}\cup\widehat{R}_{im}$ the set of roots of $\widehat{\lieg}$. For each real root $\alpha+r\delta$, we have a Lie subalgebra generated by the set $\{x^{+}_{\alpha}\otimes t^r, x^{-}_{\alpha}\otimes t^{-r}\}$ which is isomorphic to $\mathfrak{sl}_2$. 
%For $i\in I$, set $\Lambda_i=\omega_i+\langle \omega_i,\,\theta^\vee \rangle \Lambda_0 \in \widehat{\csa}^*$. 
%Set $\widehat{I}=I\cup\{0\}$, $\alpha_0=\delta-\theta$, and $\alpha_0^\vee=c-\theta^\vee$. The affine weight lattice $P$ and  the set ${\widehat{P}}^+$ of dominant integral affine weights are given by
%$$\widehat{P}=P\oplus\mathbb{Z}\Lambda_0\oplus\mathbb{Z}\delta \qquad \textup{and}\qquad {\widehat{P}}^+=\{\Lambda\in\widehat{P}:\langle\Lambda, \alpha_i^\vee \rangle{\geq 0}\,\,\forall\,\,i\in\widehat{I}\}.$$ 
The affine root lattice $\widehat{Q}$ is the $\mathbb{Z}$-span of the simple roots $\alpha_i, i\in\widehat{I}$, of $\widehat{\lieg}$, and $\widehat{Q}^+$ is defined in the obvious way. Let $\widehat{R}_{re}=\{\alpha+r\delta : \alpha \in R,\,r\in\mathbb{Z}\}$ be the set of real roots, $\widehat{R}_{im}=\{r\delta : r\in\mathbb{Z}\setminus \{0\}\}$ the set of imaginary roots and $ \widehat{R}=\widehat{R}_{re}\cup\widehat{R}_{im}$ the set of roots of $\widehat{\lieg}$. For each real root $\alpha+r\delta$, we have the Lie subalgebra of $\widehat{\lieg}$ generated by  $\{x^{+}_{\alpha}\otimes t^r, x^{-}_{\alpha}\otimes t^{-r}\}$ which is isomorphic to $\mathfrak{sl}_2$. 

Let $\afw$ be the affine Weyl group  with simple reflections $r_i (i\in\widehat{I})$. We regard $W$ naturally as a subgroup of $\afw$. 
%Let  $s_{\alpha}\in \afw$  be the reflection corresponding to $\alpha$. 
Given $\alpha \in \csa^*$, we define $t_\alpha \in GL(\csaf^*)$ by
\[
t_{\alpha}(\lambda)=\lambda+ \bform{\lambda}{\delta} \alpha -\bform{\lambda}{\alpha} \delta - \frac{1}{2} \bform{\lambda}{\delta} \bform{\alpha}{\alpha} \delta \;\;\text{ for } \lambda \in \csaf^*.
\]
The  translation subgroup $T_M$ of $\afw$ is defined by $T_M = \{t_{\alpha}\in GL(\csaf^*): \alpha \in M\}$ and we have 
\[\widehat{W} = W\ltimes T_M.\]
The extended affine Weyl group $\widetilde{W}$ is the semi-direct product 
\[\widetilde{W} =W\ltimes T_L,\]
where $T_L = \{t_{\beta}\in GL(\csaf^*): \beta \in L\}$.
We also have
$\widetilde{W} = \afw \rtimes \digautos$, where $\digautos$ is the subgroup of diagram automorphisms of $\widehat{\lieg}$. Given $w\in \widehat{W},$ let $\ell(w)$ be the length of a reduced expression of $w$. The length function $\ell$ is extended to $\widetilde{W}$ by setting 
$\ell(w\sigma)=\ell(w)$ for $w\in \afw$ and $\sigma\in\Sigma$.

The following lemma is proved in  \cite{CSVW}, and will be required  later.
\begin{lemma}\cite[Proposition 2.8]{CSVW}\label{lengthl}
Given $\lambda, \mu\in P^+$ and $w\in W$, we have 
\[\ell(t_{-\lambda} t_{-\mu} w)=\ell(t_{-\lambda})+ \ell(t_{-\mu} w).\]
\end{lemma}

For any group $G$, let $\mathbb{Z}[G]$ be the integral group ring of $G$ with basis $e^g, g\in G$. Let $I_{\delta}$ be the ideal of $\mathbb{Z}[\widehat{P}]$ obtained by setting $e^{\delta}=1.$
For a finite-dimensional semisimple $\csa$-module V, we define $\csa$-character $\textup{ch}_{\csa} V$ by
\[\textup{ch}_{\csa} V=\sum_{\mu\in\csa^*} \dim\,V_{\mu}\, e^\mu \in \mathbb{Z}[\csa^*],\]
where $V_{\mu}=\{v\in V: hv=\langle\mu,\,h\rangle \,\forall\, h\in\csa\}$. For a finite-dimensional semisimple $\widehat{\csa}$-module V, the $\widehat{\csa}$-character $\textup{ch}_{\widehat{\csa}} V\in\mathbb{Z}[\widehat{\csa}^*]$ is defined in the similar way.
\subsection{}
Let $V(\Lambda)$ be the integrable highest weight $\widehat{\lieg}$-module corresponding to a dominant integral affine weight $\Lambda$. Let $\xi_1,\ldots,\xi_p$ be a sequence of elements of $\afw(\widehat{P}^+)$. For $1\leq j\leq p$, let $\Lambda^j$ be the element of $\widehat{P}^+$ such that $\xi_j\in\afw \Lambda^j$. Define a $\widehat{\mathfrak{b}}$-submodule $D(\xi_1,\ldots,\xi_p)$ of $V(\Lambda^1)\otimes\cdots\otimes V(\Lambda^p)$ by 
\[ D(\xi_1,\ldots,\xi_p)=\mathbf{U}(\widehat{\mathfrak{b}})\big(V(\Lambda^1)_{\xi_1}\otimes\cdots\otimes V(\Lambda^p)_{\xi_p}\big).\] We call the $\widehat{\mathfrak{b}}$-module $D(\xi_1,\ldots,\xi_p)$ as a generalized Demazure module \cite{LLM, Naoi3}. When $p=1$, $D(\xi_1)$ is called as a Demazure module. We observe that  $D(\xi_1,\ldots,\xi_p)\subset D(\xi_1)\otimes\cdots\otimes D(\xi_p)$.
Note that $D(\xi_1,\ldots,\xi_p)$ is $\lieg$-stable if 
$\langle \xi_j,\,\hi\rangle \leq 0,\,\, \forall\,\, i\in I,\,\, 1\leq j\leq p.$
\subsection{}
We now recall from \cite[\S3.5]{CV} the definition of the finite-dimensional graded $\lieg[t]$-modules $D(\ell,\,\lambda)$,  $(\ell, \lambda)\in \mathbb{N}\times P^+$. 
For $\alpha\in R^{+}$ with $\langle \lambda, \alpha^{\vee}\rangle  
> 0$, let $s_\alpha, m_\alpha \in \mathbb{N}$ be the unique positive integers such that
\[ \langle \lambda, \alpha^{\vee}\rangle = (s_\alpha-1)d_\alpha \ell + m_\alpha, \quad 0<m_\alpha\leq d_\alpha \ell.\]
If $\langle \lambda, \alpha^{\vee}\rangle = 0$ we set $s_\alpha=1$ and $m_\alpha=0$. 
The module $D(\ell,\,\lambda)$ is the cyclic $\lieg[t]$-module generated by an element $w_{\ell,\,\lambda}$ with the following defining relations:
\begin{gather*}
(x^{+}_{i}\otimes t^s)\, w_{\ell,\,\lambda}=0,\quad (\hi \otimes t^s) \,w_{\ell,\,\lambda} = \langle \lambda, \hi \rangle \delta_{s,0}w_{\ell,\,\lambda}, \quad 
(x^{-}_{i} \otimes 1)^{\langle \lambda,\, \hi \rangle + 1} \,w_{\ell,\,\lambda}=0,\quad \forall \,\,s{\geq0}, i\in I,
\\
(x^{-}_{\alpha} \otimes t^{s_{\alpha}}) \,w_{\ell,\,\lambda}=0, \quad \forall\,\alpha \in R^{+},
\\
(x^{-}_{\alpha} \otimes t^{s_{\alpha}-1})^{m_{\alpha}+1} \,w_{\ell,\,\lambda}=0, \,\,\text{if}\,\, m_\alpha<d_{\alpha}\ell,\quad\forall\, \alpha\in R^{+}.
\end{gather*}
The following relations also hold in the module $D(\ell,\,\lambda)$
$$(x^{-}_{\alpha} \otimes t^{s_{\alpha}-1})^{m_{\alpha}+1} \,w_{\ell,\,\lambda}=0, \,\,\text{if}\,\, m_\alpha=d_{\alpha}\ell,\quad\forall\, \alpha\in R^{+}.$$
We declare the grade of $w_{\ell,\,\lambda}$ to be zero. Since the defining relations of $D(\ell,\,\lambda)$ are graded, it follows that $D(\ell,\,\lambda)$ is a graded $\lieg[t]$-module.
 
The following result, which gives the connection with Demazure modules,  is a combination of results \cite[Theorem 2]{CV}, \cite[Corollary 1]{FoL}, and \cite[Proposition 3.6]{Naoi1}.
\begin{proposition}\label{DvsV}\cite{CV, FoL, Naoi1}
Given $(\ell,\lambda)\in\mathbb{N}\times P^+$ and $m\in\mathbb{Z}$, let $w\in\afw, \sigma\in\Sigma$, and $\Lambda\in \widehat{P}^+$ such that 
\[w\sigma\Lambda\equiv w_0\lambda+\ell\Lambda_0+m\delta.\]
Then we have the following isomorphism of $\lieg[t]$-modules,
\[\tau_m D(\ell,\,\lambda)\cong D(w\sigma\Lambda).\] 
Under this isomorphism the generator $w_{\ell,\,\lambda}$ of $\tau_m D(\ell,\,\lambda)$ maps to a non-zero element of the weight space of $D(w\sigma\Lambda)$ of weight $w_0w\sigma\Lambda$.
\end{proposition}
\subsection{}
We recall the notion of the fusion product of finite-dimensional cyclic graded $\lieg[t]$-modules given in \cite{FL}.

Let $V$ be  a cyclic $\lieg[t]$-module generated by $v$. We define a filtration $F^{r}V, \,r\in\mathbb{Z}_{\geq 0}$ on $V$ by
\[ F^{r}V=\sum_{0\leq s \leq r}\, \mathbf{U}(\lieg[t])[s]\,v.\]
Set $F^{-1}V=\{0\}$. The associated graded space $\text{gr}\,V=\bigoplus_{r\geq 0}\, F^{r}V/ F^{r-1}V $ naturally becomes a cyclic $\lieg[t]$-module generated by the image of $v$ in $\text{gr}\,V$.

Given a $\lieg[t]$-module $V$ and a complex number $z$, we define an another  $\lieg[t]$-module action on $V$ as follows:
\[(x\otimes t^s)\,v=(x\otimes (t+z)^s)\,v, \quad x\in\lieg, \,\, v\in V, \,\,s\in \mathbb{Z}_{\geq0}.\]
Denote this new module by $V^z$. For $1\leq i\leq m$, let $V_i$ be a finite-dimensional cyclic graded $\lieg[t]$-module generated by $v_i$. Let $z_1,\ldots,z_m$ be distinct complex numbers. We denote
$\mathbf{V}={V_1}^{z_1}\otimes\cdots\otimes {V_m}^{z_m}$,
the tensor product of the corresponding  $\lieg[t]$-modules. It is known (see \cite[Proposition 1.4]{FL}) that $\mathbf{V}$ is a cyclic $\lieg[t]$-module generated by $v_1\otimes\cdots\otimes v_m$. The associated graded space $\text{gr}\mathbf{V}$ is called the fusion product of $V_1,\ldots,V_m$ w.r.t. the parameters $z_1,\ldots,z_m,$ and is denoted by ${V_1}^{z_1}*\cdots*{V_m}^{z_m}$. We denote  the image of $v_1\otimes\cdots\otimes v_m$ in $\text{gr} \mathbf{V}$ by $v_1*\cdots*v_m$. For ease of notation we shall often write $V_1*\cdots*V_m$ for ${V_1}^{z_1}*\cdots*{V_m}^{z_m}$. 
We note that $V_1*\cdots*V_m \cong_{\lieg} V_1\otimes\cdots\otimes V_m$. 
\section{The main results}
We begin this section by introducing a class of  finite-dimensional graded cyclic $\lieg[t]$-modules and then state our main results. 
%We then discuss applications of our result.

Given $k, \ell, m \in\mathbb{Z}_{\geq 0}$ such that $k\leq m$ and $\lambda\in L^+$. 
We define  $\mathbf{V}^{\ell,\,\ell\lambar}_{m,\,k\theta}$ to be the cyclic $\lieg[t]$-module 
generated by an element $v=v^{\ell,\,\ell\lambar}_{m,\,k\theta}$ with the following defining relations:
\begin{gather}
(x^{+}_{i}\otimes t^s)\,v=0,\quad (\hi\otimes t^s)\,v =\delta_{s,0} \langle \ell\lambda+k\theta, \, \hi \rangle\, v,  \quad \big(x^{-}_{i}\otimes 1\big)^{\langle \ell\lambda+k\theta, \, \hi \rangle +1} 
\,v=0,\quad\forall \,\, s\geq0, i\in I,\label{mr:dfr1}
%\\
%\big(x^{-}_{\alpha}\otimes 1\big)^{\langle \ell\lambda+k\theta, \, \alpha^{\vee} \rangle +1} \,v=0, \quad \forall \,\, \alpha\in R^+,\label{mr:dfr2}
\\
\big(x^{-}_{\alpha}\otimes t^{(\lambda+\theta|\alpha)}\big)\,v=0, \quad \forall \,\, \alpha\in R^+,\label{mr:dfr3}
\\
\big(x^{-}_{\alpha}\otimes t^{(\lambda|\alpha)}\big)^{\langle (k-i)\theta,\, \alpha^{\vee} \rangle +1} \big(x^{-}_{\theta}\otimes t^{(\lambda|\theta)+1}\big)^i\,v=0, \quad \forall \,\,\alpha\in R^+,\,\,\, 0\leq i \leq k,\label{redrel}
\\
\big(x_{\theta}^-\otimes t^{(\lambda|\theta)+1}\big)^{2k-m+1}\,v=0, \quad \text{if}\,
\,\, m\leq 2k, \label{mr:dfr5}
\\
\big(x_{\theta}^-\otimes t^{(\lambda|\theta)+1}\big)\,v=0, \quad \text{if} \,\,\, m\geq 2k.\label{mr:dfr6}
\end{gather}
The relations \eqref{mr:dfr1} guarantee that the module $\mathbf{V}^{\ell,\,\ell\lambar}_{m,\,k\theta}$ is finite-dimensional (cf. \cite{CP}).
In particular this gives,
\beq\label{mr:dfr7}
\big(x^{-}_{\alpha}\otimes 1\big)^{\langle \ell\lambda+k\theta, \, \alpha^{\vee} \rangle +1} \,v=0, \quad \forall \,\, \alpha\in R^+.
\eeq
We declare the grade of $v$ to be zero. Since the defining relations of $\mathbf{V}^{\ell,\,\ell\lambar}_{m,\,k\theta}$ are graded, 
it follows that $\mathbf{V}^{\ell,\,\ell\lambar}_{m,\,k\theta}$  is a graded $\lieg[t]$-module. 
\begin{remark}
We prove in \S\ref{connectiontoCV} that the relations in \eqref{redrel}, when $\alpha\in R^+\setminus\{\theta\}$ and $1\leq i\leq k$, are redundant (see Lemma \ref{lemCV}).
\end{remark}

We observe that
\beq\label{obs0}
\mathbf{V}^{\ell,\,\ell\lambar}_{m,\,k\theta}\cong_{\lieg[t]} \mathbf{V}^{\ell,\,\ell\lambar}_{2k,\,k\theta},\quad\forall \,\,\,m\geq 2k,
\eeq
and
\beq\label{obs1}
\mathbf{V}^{\ell,\,\ell\lambar}_{0,\,0\theta}\cong_{\lieg[t]}D(\ell,\,\ell\lambda).
\eeq

The following theorems are the main results of this paper.
\begin{theorem}\label{MTone}
Let $k, \ell, m\in\mathbb{N}$ and $\lambda\in L^+$.
\be
\item\label{pone}
If $\ell\geq 2k$,  then there exists a short exact sequence of $\lieg[t]$-modules,
\[
0\rightarrow \tau_{(\lambda|\theta)+1} {\mathbf{V}}^{\ell,\,\ell\lambar}_{k-1,\,(k-1)\theta}
\xrightarrow{{\phi^{-}_1}}{\mathbf{V}}^{\ell,\,\ell\lambar}_{k,\,k\theta}
\xrightarrow{{\phi^{+}_k}} {\mathbf{V}}^{\ell,\,\ell\lambar}_{\ell,\,k\theta}\rightarrow 0.
\]
\item\label{ptwo}
If  $\ell\geq 2k$ and $k<m\leq 2k$, then there exists a short exact sequence of $\lieg[t]$-modules,
\[0\rightarrow \tau_{(2k-m+1)((\lambda|\theta)+1)} {\mathbf{V}}^{\ell,\,\ell\lambar}_{m-k-1,\,(m-k-1)\theta}
\xrightarrow{{\psi}^{-}}{\mathbf{V}}^{\ell,\,\ell\lambar}_{k,\,k\theta}
\xrightarrow{\psi^+} {\mathbf{V}}^{\ell,\,\ell\lambar}_{m,\,k\theta}\rightarrow 0.\]
\item\label{pthree}
If $\ell\leq 2k$ and $k\leq m<\ell$, then  there exists a short exact sequence of $\lieg[t]$-modules,
\begin{gather*}
0\rightarrow \tau_{(2k-\ell+1)((\lambda|\theta)+1)} {\mathbf{V}}^{\ell,\,\ell\lambar}_{\ell+m-2k-1,\,(\ell-k-1)\theta}
\xrightarrow{{\phi^{-}_2}}{\mathbf{V}}^{\ell,\,\ell\lambar}_{m,\,k\theta}
\xrightarrow{{\phi^{+}_m}} {\mathbf{V}}^{\ell,\,\ell\lambar}_{\ell,\,k\theta}\rightarrow 0.
\end{gather*}
\ee
\end{theorem}
\begin{theorem}\label{MTtwo}Let $k, \ell, m\in\mathbb{N}$ be such that  $\ell\geq m\geq k$. 
Let $\lambda_1,\ldots,\lambda_p$ be a sequence of elements of ${L^+}$, and denote $\lambda=\lambda_1+\cdots+\lambda_p$.
Then we have the following isomorphisms of $
\lieg[t]$-modules:
\begin{gather*}
\begin{split}
&D(\ell,\,\ell\lambda_1)*\cdots*D(\ell,\,\ell\lambda_p)*D(m,\,k\theta)\\&\cong{\mathbf{V}}^{\ell,\,\ell\lambar}_{m,\,k\theta}
\cong \begin{cases} D\big(t_{w_0\lambda} (\ell-m)\Lambda_0,\, t_{w_0(\lambda+\theta)} (m\Lambda_0+(m-k)\theta)\big), &   m\leq 2k,\\
D\big(t_{w_0\lambda}(\ell-m)\Lambda_0,\, t_{w_0\lambda}w_0 (m\Lambda_0+k\theta)\big), &  m\geq 2k.\end{cases}
\end{split}\end{gather*}
\end{theorem}

Theorems \ref{MTone} and \ref{MTtwo} are proved in \S\ref{s:proof}.

%\subsection{} In this subsection we discuss the applications of Theorem \ref{MT1}.

%The following corollary gives  short exact sequences of fusion products of Demazure modules and of generalized Demazure modules respectively. 
%Their proves are immediate from Theorem \ref{MT1}.
The following corollary is immediate from Theorems \ref{MTone} and \ref{MTtwo}.
\begin{corollary}\label{exactfg}
Let $k, \ell, m\in\mathbb{N}$ and $\lambda_1,\ldots,\lambda_p\in {L^+}$. Denote $\lambda=\lambda_1+\cdots+\lambda_p$. 
\be
\item
If $\ell\geq 2k$, then there exist short exact sequences of $\lieg[t]$-modules,
\begin{align*}
0&\rightarrow \tau_{(\lambda|\theta)+1} \big(D(\ell,\,\ell\lambda_1)*\cdots*D(\ell,\,\ell\lambda_p)*D(k-1,\,(k-1)\theta)\big)
\\
&\rightarrow D(\ell,\,\ell\lambda_1)*\cdots*D(\ell,\,\ell\lambda_p)*D(k,\,k\theta)
\\
&\rightarrow 
 D(\ell,\,\ell\lambda_1)*\cdots*D(\ell,\,\ell\lambda_p)*D(\ell,\,k\theta) \rightarrow 0\end{align*}
 and
 \begin{align*}
0&\rightarrow \tau_{(\lambda|\theta)+1} \,D\big(t_{w_0\lambda}(\ell-k+1)\Lambda_0,\, t_{w_0(\lambda+\theta)}(k-1)\Lambda_0\big)
\\
&\rightarrow D\big(t_{w_0\lambda}(\ell-k)\Lambda_0,\, t_{w_0(\lambda+\theta)}k\Lambda_0\big)
\\
&\rightarrow 
 D\big(t_{w_0(\lambda+\theta)}(\ell\Lambda_0+(\ell-k)\theta)\big) \rightarrow 0.
\end{align*}
\item
If $\ell\geq2k$ and $k< m\leq 2k$, then there exist short exact sequences of $\lieg[t]$-modules,
\begin{align*}
0&\rightarrow \tau_{(2k-m+1)((\lambda|\theta)+1)} \big(D(\ell,\,\ell\lambda_1)*\cdots*D(\ell,\,\ell\lambda_p)*D(m-k-1,\,(m-k-1)\theta)\big)
\\
&\rightarrow D(\ell,\,\ell\lambda_1)*\cdots*D(\ell,\,\ell\lambda_p)*D(k,\,k\theta)
\\
&\rightarrow 
 D(\ell,\,\ell\lambda_1)*\cdots*D(\ell,\,\ell\lambda_p)*D(m,\,k\theta) \rightarrow 0
 \end{align*}
 and
 \begin{align*}
0&\rightarrow \tau_{(2k-m+1)((\lambda|\theta)+1)} \,D\big(t_{w_0\lambda}(\ell-m+k+1)\Lambda_0,\, t_{w_0(\lambda+\theta)}(m-k-1)\Lambda_0\big)
\\
&\rightarrow D\big(t_{w_0\lambda}(\ell-k)\Lambda_0,\, t_{w_0(\lambda+\theta)}k\Lambda_0\big)
\\
&\rightarrow 
 D\big(t_{w_0\lambda}(\ell-m)\Lambda_0,\, t_{w_0(\lambda+\theta)}(m\Lambda_0+(m-k)\theta)\big)\rightarrow 0.
\end{align*}
\item
If $\ell\leq2k$ and $k\leq m<\ell$, then there exist  short exact sequences of $\lieg[t]$-modules,
\begin{align*}
0&\rightarrow \tau_{(2k-\ell+1)((\lambda|\theta)+1)} \big(D(\ell,\,\ell\lambda_1)*\cdots*D(\ell,\,\ell\lambda_p)*D(\ell+m-2k-1,\,(\ell-k-1)\theta)\big)
\\
&\rightarrow D(\ell,\,\ell\lambda_1)*\cdots*D(\ell,\,\ell\lambda_p)*D(m,\,k\theta)
\\
&\rightarrow 
 D(\ell,\,\ell\lambda_1)*\cdots*D(\ell,\,\ell\lambda_p)*D(\ell,\,k\theta) \rightarrow 0
 \end{align*}
 and
 \begin{align*}
0&\rightarrow \tau_{(2k-\ell+1)((\lambda|\theta)+1)} \,D\big(t_{w_0\lambda}(2k-m+1)\Lambda_0,\, t_{w_0(\lambda+\theta)}((\ell+m-2k-1)\Lambda_0+(m-k)\theta)\big)
\\
&\rightarrow
D\big(t_{w_0\lambda}(\ell-m)\Lambda_0,\, t_{w_0(\lambda+\theta)}(m\Lambda_0+(m-k)\theta)\big)
\\
&\rightarrow 
 D\big(t_{w_0(\lambda+\theta)}(\ell\Lambda_0+(\ell-k)\theta)\big) \rightarrow 0.
\end{align*}
\ee
\end{corollary}
%\begin{proof}
%The proof is immediate from Theorem \ref{MT1}.
%\end{proof}
The following corollary gives a surjective morphism between two fusion products of Demazure modules.
\begin{corollary}\label{Schur}
Let $(\ell\geq m\geq k)$ and $(\ell^\prime\geq m^\prime\geq k^\prime)$ be two sequences of positive integers with $k^\prime\leq k$. Let
$\lambda_1,\ldots,\lambda_p$ and $\lambda^\prime_1,\ldots,\lambda^\prime_p$ be two sequences of elements in $L^+$. Denote $\lambda=\sum_{i=1}^{p}\lambda_i$ and $\lambda^\prime=\sum_{i=1}^{p}\lambda^\prime_i$. Suppose that
\begin{enumerate}
\item[(a)]\label{hypoa}
$\ell\lambda+k\theta=\ell^\prime\lambda^\prime+k^\prime\theta,$
\item[(b)]\label{hypob}
$(\lambda|\alpha)\geq(\lambda^\prime|\alpha), \quad\forall\,\,\alpha\in R^+.$
\end{enumerate}
If $(\lambda|\theta)=(\lambda^\prime|\theta)$, then we further assume the following holds;
\begin{enumerate}
\item[(c)]\label{hypoc} if $m\geq 2k$, then  $m^\prime\geq 2k^\prime$,
\item[(d)]\label{hypod} if $m\leq 2k$ and $m^\prime\leq 2k^\prime$, then $2k-m\geq2k^\prime-m^\prime$.
\end{enumerate}
%Suppose that $k\leq k^\prime$.
Then  there exist surjective morphisms of $\lieg[t]$-modules,
\beq\label{Schurp1}
 {\mathbf{V}}^{\ell,\,\ell\lambar}_{m,\,k\theta}\twoheadrightarrow{\mathbf{V}}^{\ell^\prime,\,\ell^\prime\lamprbar}_{m^\prime,\,k^\prime\theta}\eeq and
\beq\label{Schurp2}
D(\ell,\,\ell\lambda_1)*\cdots*D(\ell,\,\ell\lambda_p)*D(m,\,k\theta)\twoheadrightarrow D(\ell^\prime,\,\ell^\prime\lambda^\prime_1)*\cdots*D(\ell^\prime,\,\ell^\prime\lambda^\prime_p)*D(m^\prime,\,k^\prime\theta).
\eeq
\end{corollary}
\begin{proof}
To prove  \eqref{Schurp1}, we  show that the generator $v^\prime=v^{\ell^\prime,\,\ell^\prime\lamprbar}_{m^\prime,\,k^\prime\theta}$ of ${\mathbf{V}}^{\ell^\prime,\,\ell^\prime\lamprbar}_{m^\prime,\,k^\prime\theta}$ satisfies the defining relations of 
 ${\mathbf{V}}^{\ell,\,\ell\lambar}_{m,\,k\theta}$. The relations \eqref{mr:dfr1} are clear from hypothesis (a). 
 It is easy to see that the relations \eqref{mr:dfr3} follow by using hypothesis (b). 
 We  prove the relations \eqref{redrel} by considering two cases. First suppose that $(\lambda|\alpha)>(\lambda^\prime|\alpha)$. Then, they follow from the hypothesis  $k^\prime\leq k$, by using the relations
 $$\big(x^{-}_{\alpha}\otimes t^{(\lambda^\prime|\alpha)+1}\big)\,v^\prime=0, \quad\textup{if}\,\,\alpha\neq\theta,\qquad\quad
 \big(x^{-}_{\theta}\otimes t^{(\lambda^\prime|\theta)+1}\big)^{k^\prime+1}\,v^\prime=0,$$
 and
 \beq\label{pfsec3e1}
 \big(x^{-}_{\theta}\otimes t^{(\lambda^\prime|\theta)+2}\big)\,v^\prime=0.
 \eeq
 In the second case, using hypothesis (b), we have  $(\lambda|\alpha)=(\lambda^\prime|\alpha)$. Then, they follow from the hypothesis  $k^\prime\leq k$ and $(\lambda|\theta)\geq(\lambda^\prime|\theta)$, by using \eqref{pfsec3e1}.
 
 We prove the remaining relations \eqref{mr:dfr5} and \eqref{mr:dfr6} by considering two cases. 
 First suppose that
 $(\lambda|\theta)>(\lambda^\prime|\theta)$. Then, they are immediate from  \eqref{pfsec3e1}. 
 %$\big(x^{-}_{\theta}\otimes t^{(\lambda^\prime|\theta)+2}\big)\,v^\prime=0$.
 In the second case, we have $(\lambda|\theta)=(\lambda^\prime|\theta)$. 
 Then, the relation \eqref{mr:dfr5} (resp. \eqref{mr:dfr6}) in this case is
 clear from hypothesis (d) (resp. hypothesis (c)). Hence  \eqref{Schurp1}. 
 The proof of  \eqref{Schurp2} is immediate from  \eqref{Schurp1} by using Theorem \ref{MTtwo}.
\end{proof}
The next corollary, which gives a surjective morphism between two tensor products of Demazure modules, may be viewed as a generalization of the Schur positivity \cite{CFS}.
\begin{corollary}\label{Schurpm}
Let $(\ell_1\geq \ell_2\geq\cdots\geq \ell_{p}\geq0)$ and $(m_1\geq m_2\geq\cdots\geq m_{p}\geq0)$ be two partitions of a positive integer. Suppose that  $\ell_i+\cdots+\ell_p\geq m_i+\cdots+m_p$  for each $1\leq i\leq p$. Then 
%such that $\bm{m}\preceq{\bm\ell}$, where $\preceq$
%is the reverse dominance order on partitions of $N$, 
there exists a surjective morphism of $\lieg$-modules,
$$D(\ell_1,\,\ell_1\theta)\otimes D(\ell_2,\,\ell_2\theta)\otimes\cdots\otimes D(\ell_{p},\,\ell_{p}\theta)\twoheadrightarrow D(m_1,\,m_1\theta)\otimes D(m_2,\,m_2\theta)\otimes\cdots\otimes D(m_{p},\,m_{p}\theta).$$
\end{corollary}
\begin{proof}
It is easy to see by using similar arguments of the proof of \cite[Theorem 1 (ii)]{CFS} that it is enough to prove it for $p=2$. The proof in $p=2$ case follows from 
%the $p=1$ case of 
Corollary \ref{Schur}.
\end{proof}
\section{Proof of the main results}\label{s:proof}
The main goal of this section is to prove Theorems \ref{MTone} and \ref{MTtwo}. The proof uses a result from \cite{CSVW} and the character of generalized Demazure modules.
%In the course of the proof, we discuss about the Demazure modules, Demazure operators,  Demazure character formula, and  generalized Demazure modules. %At the end we prove Theorem \ref{MT1}.
\subsection{}
In this subsection, we give the $\lieg$-module decomposition for the modules $D(\ell,\,k\theta)$, $\ell\geq k$.
%The following proposition gives the $\lieg$-module decomposition for the modules $D(\ell, k\theta), k\leq \ell\leq 2k$.
\begin{proposition}\label{DlkTP}
Let $k,\ell\in\mathbb{N}$ be such that $k\leq \ell\leq 2k$. For $i\in\mathbb{Z}_{\geq 0}$, the subspace of $D(\ell,\,k\theta)$ of grade $i$ is given by
\begin{gather*}
D(\ell,k\theta)[i]=\begin{cases} \mathbf{U}(\lieg)(x^{-}_{\theta}\otimes t)^{i}\,w_{\ell,\,k\theta}\cong_{\lieg} V((k-i)\theta),  &  0\leq i \leq 2k-\ell,\\
\{0\},  &  i>2k-\ell.\end{cases}
\end{gather*}
In particular,
\[ D(\ell, k\theta)\cong_{\lieg} V(k\theta)\oplus V((k-1)\theta)\oplus\cdots\oplus V((\ell-k)\theta).\]
\end{proposition}
\begin{proof}
We prove this by considering the Demazure module corresponding to the module $D(\ell,k\theta)$. Observe that the condition $k\leq \ell\leq 2k$ is equivalent to 
$\ell\Lambda_0+(\ell-k)\theta\in \widehat{P}^{+}$, and 
\[t_{w_0(\theta)}(\ell\Lambda_0+(\ell-k)\theta)=\ell\Lambda_0+w_0(k\theta)-
(2k-\ell)\delta.\] Hence by Proposition \ref{DvsV}, we get
\beq\label{DisoV}
D(\ell,k\theta)\cong_{\lieg[t]}D(t_{w_0(\theta)}(\ell\Lambda_0+(\ell-k)\theta)).
\eeq
Under this isomorphism, the generator $w_{\ell,\,k\theta}$ of $D(\ell,\,k\theta)$ maps to a non-zero element $v$ of the weight space of $D(t_{w_0(\theta)}(\ell\Lambda_0+(\ell-k)\theta))$ of weight ${\ell\Lambda_0+k\theta-(2k-\ell)\delta}.$ 
Considering the $\mathfrak{sl}_2$ copy associated to the real root $\theta-\delta$, it follows from the standard $\mathfrak{sl}_2$ arguments that
\[(x^{-}_{\theta}\otimes t)^{2k-\ell}\,v\neq 0,\]
since $(\ell\Lambda_0+k\theta-(2k-\ell)\delta| \theta-\delta)=2k-\ell\geq 0.$ In particular,
\[(x^{-}_{\theta}\otimes t)^{2k-\ell}\,w_{\ell,\,k\theta}\neq 0.\]
The proof now follows by using the following relations which hold in the module $D(\ell,k\theta)$
\[(x^{-}_{\alpha}\otimes t)\,w_{\ell,\,k\theta}=0, \,\, \forall\,\,\alpha\in R^+\setminus \{\theta\},\qquad (x^{-}_{\theta}\otimes t^2)\,w_{\ell,\,k\theta}=0, \qquad (x^{-}_{\theta}\otimes t)^{2k-\ell+1}\,w_{\ell,\,k\theta}=0,
\]
and
\[\mathfrak{n}^+(x^{-}_{\theta}\otimes t)^{i}\,w_{\ell,\,k\theta}=0, \quad\forall\,\,\, i\in\mathbb{Z}_{\geq 0}.\]
%This completes the proof.
\end{proof}
We record below an easy fact, for later use. 
\beq\label{lgeq2k}
D(\ell,k\theta)\cong_{\lieg[t]}D(w_0(\ell\Lambda_0+k\theta))\cong_{\lieg[t]} \textup{ev}_0\,V(k\theta), \quad\forall\,\ell\geq 2k.
\eeq
\subsection{} We now recall the result \cite[Theorem 1]{CSVW} for $\lambda^0=k\theta$. Its proof uses the conclusion of \cite[Proposition 3.5]{CSVW} when $\lambda=\lambda^0$,
due to which it has some constrains in its hypothesis (see \cite[Remark~3.4]{CSVW}).
 We observe that 
  $$t_{\theta}(\ell\Lambda_0+w_0k\theta)=\ell\Lambda_0+(\ell-k)\theta+(2k-\ell)\delta\in\widehat{P}^+, \quad\textup{if}\,\,k\leq\ell\leq 2k,$$ and 
  $$w_0(\ell\Lambda_0+w_0k\theta)\in \widehat{P}^+, \quad\textup{if}\,\,\ell\geq 2k.$$ Hence the conclusion of \cite[Proposition 3.5]{CSVW}
  when $\lambda=k\theta$ is satisfied, so we can remove the constrains in the hypothesis of \cite[Theorem 1]{CSVW} when $\lambda^0=k\theta$.
  %Now the following result is a consequence of \cite[Theorem 1 and Remark 3.4]{CSVW}.
%The following result is a consequence of \cite[Theorem 1 and Remark 3.5]{CSVW}, by using \eqref{DisoV} and \eqref{lgeq2k}.
\begin{theorem}\cite{CSVW}\label{CSVWT}
Let $k\in\mathbb{Z}_{\geq0}$ and  $\ell\in\mathbb{N}$ be such that $\ell\geq k$. Let $\lambda_1,\ldots,\lambda_p$ be a sequence of elements of $L^+$. Then we have the following isomorphism of $\lieg[t]$-modules,
\[D(\ell,\,\ell\lambda_1)*\cdots*D(\ell,\,\ell\lambda_p)*D(\ell,\,k\theta)\cong D(\ell,\,\ell(\lambda_1+\cdots+\lambda_p)+k\theta).\] 
Under this isomorphism, the generator $w_{\ell,\,\ell\lambda_1}*\cdots*w_{\ell,\,\ell\lambda_p}*w_{\ell,\,k\theta}$ maps to the generator\\
$w_{\ell,\,\ell(\lambda_1+\cdots+\lambda_p)+k\theta}$.
\end{theorem}
\subsection{} In this subsection, we prove that for $\ell\geq k$, the module $\mathbf{V}^{\ell,\,\ell\lambar}_{\ell,\,k\theta}$ is 
isomorphic to the Demazure module $D(\ell,\,\ell\lambda+k\theta)$.
%The following lemma will be useful.
\begin{lemma}\label{l2}
Let $k, \ell, m\in\mathbb{N}$ be such that $\ell\geq m\geq k$ and $\lambda_1,\ldots,\lambda_p\in {L^+}$. 
Denote $\lambda=\lambda_1+\cdots+\lambda_p$ and $w= w_{\ell,\,\ell\lambda_1}*\cdots*w_{\ell,\,\ell\lambda_p}*w_{m,\,k\theta}$. 
Then the following relations hold in the module $D(\ell,\,\ell\lambda_1)*\cdots*D(\ell,\,\ell\lambda_p)*D(m,\,k\theta)$
\be
\item\label{lemp1}
$\big(x_{\alpha}^-\otimes t^{(\lambda+\theta|\alpha)}\big)\,w=0,\quad \forall \,\, \alpha\in R^+,$
\item\label{lemp2}
$\big(x^{-}_{\alpha}\otimes t^{(\lambda|\alpha)}\big)^{\langle (k-i)\theta,\, \alpha^{\vee} \rangle +1} \big(x^{-}_{\theta}\otimes t^{(\lambda|\theta)+1}\big)^i\,w=0, \quad \forall \,\,\alpha\in R^+,\,\,\, 0\leq i \leq k,$
\item\label{lemp3}
$\big(x_{\theta}^-\otimes t^{(\lambda|\theta)+1}\big)^{2k-m+1}\,w=0, \quad \text{if}\,
\,\, m\leq 2k,$
\item\label{lemp4}
$\big(x_{\theta}^-\otimes t^{(\lambda|\theta)+1}\big)\,w=0, \quad \text{if} \,\,\, m\geq 2k.$
\ee
\end{lemma}
\begin{proof}
Let $z_1,\ldots,z_p,z_{p+1}$ be the distinct complex numbers which define the fusion product. In the corresponding tensor product, we have
\begin{equation*}
\begin{split}
&\big(x_{\alpha}^-\otimes (t-z_1)^{(\lambda_1|\alpha)}\cdots(t-z_p)^{(\lambda_p|\alpha)}(t-z_{p+1})^{(\theta|\alpha)}\big)\,(w_{\ell,\,\ell\lambda_1}\otimes\cdots\otimes w_{\ell,\,\ell\lambda_p}\otimes w_{m,\,k\theta})\\
&=\sum_{j=1}^{p+1}\big(w_{\ell,\,\ell\lambda_1}\otimes\cdots\otimes(x_{\alpha}^-\otimes t^{(\lambda_j|\alpha)}f_j(t)\,w_{\ell,\,\ell\lambda_j})\otimes\cdots\,w_{\ell,\,\ell\lambda_p}\otimes w_{m,\,k\theta}\big)=0,
\end{split}
\end{equation*}
where  $f_j(t)=\prod_{i\neq j}(t+z_j-z_i)^{(\lambda_i|\alpha)}$, $\lambda_{p+1}=\theta$, and the last equality follows by using the relations
$$(x_{\alpha}^-\otimes t^{(\lambda_j|\alpha)})\,w_{\ell,\,\ell\lambda_j}=0, \,\,\forall\,\,1\leq j\leq p \qquad \textup{and} \qquad (x_{\alpha}^-\otimes t^{(\theta|\alpha)})\,w_{m,\,k\theta}=0.$$
Now part \eqref{lemp1} is immediate. The proof of part \eqref{lemp2} (resp. part \eqref{lemp3}, part \eqref{lemp4})  
is identical by using the relation  
$$(x^{-}_{\alpha}\otimes 1)^{\langle (k-i)\theta,\, \alpha^{\vee} \rangle +1} (x^{-}_{\theta}\otimes t)^i\,w_{m,\,k\theta}=0 \quad (\textup{follows from Proposition \ref{DlkTP}})$$
$$(\textup{resp.}\,\, (x_{\theta}^-\otimes t)^{2k-m+1}\,w_{m,\,k\theta}=0 \quad(\textup{since}\,\, m\leq 2k),\qquad (x_{\theta}^-\otimes t)\,w_{m,\,k\theta}=0 \quad(\textup{since}\,\, m\geq 2k)),$$
and we omit the details.
\end{proof}
The following proposition gives explicit defining relations for the modules $D(\ell,\ell\lambda+k\theta)$, $\ell\geq k$.
\begin{proposition}\label{Dllamkt}
Let $k,\ell\in\mathbb{N}$ be such that $\ell\geq k$ and  $\lambda\in L^+$. Let $w={w_{\ell,\,\ell\lambda+k\theta}}$ be the generator of the module $D(\ell,\ell\lambda+k\theta)$. Then:
\be
\item The following are the defining relations for the module $D(\ell,\ell\lambda+k\theta)$
\begin{gather*}
%\mathfrak{n}^{+}[t]\, w=0,\qquad (h \otimes t^s) \,w = \delta_{s,0} \langle \ell\lambda+k\theta,\, h \rangle  w, \quad \forall \,s\geq0,\, h\in\csa,
%\\
(x^{+}_{i}\otimes t^s)\,w=0,\,\,\,\, (\hi\otimes t^s)\,w =\delta_{s,0} \langle \ell\lambda+k\theta, \, \hi \rangle\, w,  \,\,\,\, \big(x^{-}_{i}\otimes 1\big)^{\langle \ell\lambda+k\theta, \, \hi \rangle +1} 
\,w=0, \quad\forall \,\, s\geq0, i\in I,\\
\big(x_{\alpha}^{-}\otimes t^{(\lambda+\theta|\alpha)}\big)\,w=0, \quad\forall\,\alpha \in R^{+},
\\
\big(x_{\alpha}^{-}\otimes t^{(\lambda|\alpha)}\big)^{\langle k\theta,\,\alpha^{\vee}\rangle+1}\,w=0, \quad \forall\,\alpha \in R^{+},
\\
{\big(x^{-}_{\theta} \otimes t^{(\lambda|\theta)+1}\big)}^{2k-\ell+1} \,w=0,\quad \textup{if}\,\,\ell\leq 2k,
\\
\big(x^{-}_{\theta} \otimes t^{(\lambda|\theta)+1}\big) \,w=0, \quad \textup{if}\,\,\ell\geq 2k.
\end{gather*}
\item
The following relations also hold in the module $D(\ell,\ell\lambda+k\theta)$
\[\big(x^{-}_{\alpha}\otimes t^{(\lambda|\alpha)}\big)^{\langle (k-i)\theta,\, \alpha^{\vee} \rangle +1} \big(x^{-}_{\theta}\otimes t^{(\lambda|\theta)+1}\big)^i\,w=0, 
\quad \forall \,\,\alpha\in R^+,\,\, 1\leq i \leq k.
\]
\ee
\end{proposition}
\begin{proof}
For $\alpha\in R^+$ with $\langle \ell\lambda+k\theta, \alpha^{\vee}\rangle>0$, we first write down the positive integers $s_\alpha, m_\alpha$ such that
\[ \langle \ell\lambda+k\theta, \alpha^{\vee}\rangle = (s_\alpha-1)d_\alpha \ell + m_\alpha, \quad 0<m_\alpha\leq d_\alpha \ell.\]
If $(\theta|\alpha)=0$ (resp. $(\theta|\alpha)=1$), then $s_\alpha=(\lambda|\alpha)$ (resp. $s_\alpha=(\lambda|\alpha)+1$) and $m_\alpha=d_\alpha\ell$ (resp. $m_\alpha=d_\alpha k$). 
We now consider the remaining case when $\alpha=\theta$. If $\ell<2k$ (resp. $\ell\geq2k$), then $s_{\theta}=(\lambda|\theta)+2$ (resp. $s_\theta=(\lambda|\theta)+1$) and 
$m_\theta=2k-\ell$ (resp. $m_\theta=2k$).
The proof now follows from the definition of the module $D(\ell,\, \ell\lambda+k\theta)$, by using  Theorem \ref{CSVWT} and Lemma \ref{l2}.
\end{proof}
%The following proposition is useful later.
\begin{proposition}\label{Prop}
Given $k,\ell\in\mathbb{N}$ such that $\ell\geq k$ and $\lambda\in L^+$, we have the following isomorphism of $\lieg[t]$-modules,
\[
\mathbf{V}^{\ell,\,\ell\lambar}_{\ell,\,k\theta}\cong D(\ell,\ell\lambda+k\theta).\]
\end{proposition}
\begin{proof}
The proof is immediate from Proposition \ref{Dllamkt}.
\end{proof}
\subsection{}
We now establish the existence of the maps  ${\phi^{+}_m}$, ${\phi^{-}_1}$ and ${\phi^{-}_2}$ from Theorem \ref{MTone}.

The following proposition, which gives the existence of ${\phi^{+}_m}$, is trivially checked.
\begin{proposition}\label{existp}
Let $k, \ell, m\in\mathbb{N}$ be such that $\ell\geq m\geq k$ and $\lambda\in L^+$. 
%Let  $\lambar=(\lambda_1,\ldots,\lambda_p)\in ({L^+})^{p}$, for $p\in\mathbb{N},$ and denote that $\lambda=\lambda_1+\cdots+\lambda_p$. 
Then the map 
\[{\phi^{+}_m}:{\mathbf{V}}^{\ell,\,\ell\lambar}_{m,\,k\theta}
\rightarrow {\mathbf{V}}^{\ell,\,\ell\lambar}_{\ell,\,k\theta}\quad \text{such that} \quad {\phi^{+}_m}(v^{\ell,\,\ell\lambar}_{m,\,k\theta})= v^{\ell,\,\ell\lambar}_{\ell,\,k\theta},\]
is a surjective morphism of $\lieg[t]$-modules with 
\[ \textup{ker}\, {\phi^{+}_m}= \begin{cases} \mathbf{U}(\lieg[t])\,\big(x_\theta^{-}\otimes t^{(\lambda|\theta)+1}\big)^{2k-\ell+1}\,v^{\ell,\,\ell\lambar}_{m,\,k\theta},
&   \ell\leq 2k, \\ \mathbf{U}(\lieg[t])\,\big(x^-_{\theta}\otimes t^{(\lambda|\theta)+1}\big)\,v^{\ell,\,\ell\lambar}_{m,\,k\theta}, & \ell\geq 2k.
\end{cases}
\]
\end{proposition}
%The following lemma helps in proving the existence of ${\phi^{-}_1}$ and ${\phi^{-}_2}$.
\begin{lemma}\label{l1}
Given $\lambda\in L^+$, 
%$\lambar=(\lambda_1,\ldots,\lambda_p)\in ({L^+})^{p}$, for $p\in\mathbb{N}$, and denote that $\lambda=\lambda_1+\cdots+\lambda_p.$ 
we have the following:
\be
\item\label{lemmapart1}
For $k,\ell,m\in\mathbb{N}$ such that $k\leq m<\ell\leq 2k$,  the following relations hold in the module ${\mathbf{V}}^{\ell,\,\ell\lambar}_{m,\,k\theta}.$
\be
\item\label{lemmapart1a}
$(x_{\alpha}^+\otimes t^s)(x_{\theta}^-\otimes t^{(\lambda|\theta)+1})^{2k-\ell+1}\,v^{\ell,\,\ell\lambar}_{m,\,k\theta}=0,\quad \forall \,\,s\geq 0,\,\, \alpha\in R^+,$
\item\label{lemmapart1b}
$(x_{\alpha}^-\otimes 1)^{\langle \ell\lambda+(\ell-k-1)\theta,\, \alpha^\vee\rangle+1}(x_{\theta}^-\otimes t^{(\lambda|\theta)+1})^{2k-\ell+1}\,v^{\ell,\,\ell\lambar}_{m,\,k\theta}=0, \quad \forall \,\, \alpha\in R^+.$
\ee
\item\label{lemmapart2}
For  $k, \ell\in\mathbb{N}$ such that $\ell\geq k$,  the following relations hold in the module ${\mathbf{V}}^{\ell,\,\ell\lambar}_{k,\,k\theta}.$
\be
\item
$(x_{\alpha}^+\otimes t^s)(x_{\theta}^-\otimes t^{(\lambda|\theta)+1})\,v^{\ell,\,\ell\lambar}_{k,\,k\theta}=0,\quad \forall \,\,s\geq 0,\,\, \alpha\in R^+,$
\item
$(x_{\alpha}^-\otimes 1)^{\langle \ell\lambda+(k-1)\theta,\, \alpha^\vee\rangle+1}(x_{\theta}^-\otimes t^{(\lambda|\theta)+1})\,v^{\ell,\,\ell\lambar}_{k,\,k\theta}=0, \quad \forall\,\, \alpha\in R^+.$
\ee
\ee
\end{lemma}
\begin{proof}
Set $v=v^{\ell,\,\ell\lambar}_{m,\,k\theta}$.
To prove part \eqref{lemmapart1a}, since $(x_{\alpha}^+\otimes t^s)\,v=0$, it is enough to prove that
$$[x_{\alpha}^+\otimes t^s, (x_{\theta}^-\otimes t^{(\lambda|\theta)+1})^{2k-\ell+1}]\,v=0.$$ This means to show that
\beq\label{pflemeq1}
\sum_{j=0}^{2k-\ell} \big(x_{\theta}^-\otimes t^{(\lambda|\theta)+1}\big)^j 
\big([x_{\alpha}^+, x_{\theta}^-]\otimes t^{(\lambda|\theta)+s+1}\big) \big(x_{\theta}^-\otimes t^{(\lambda|\theta)+1}\big)^{2k-\ell-j}\,v=0.
\eeq
If $(\theta|\alpha)=0$, then $[x_{\alpha}^+, x_{\theta}^-]=0$ and hence \eqref{pflemeq1}. 
If $(\theta|\alpha)=1$, then $\theta-\alpha\in R^+$ and $[x_{\alpha}^+, x_{\theta}^-]$ is 
a non-zero scalar multiple of $x_{\theta-\alpha}^-$. The proof of part \eqref{pflemeq1} in this case follows now from the relation
$(x_{\theta-\alpha}^-\otimes t^{(\lambda+\theta|\theta-\alpha)})\,v=0$.
In the remaining case when $\alpha=\theta$, it follows from the relations 
$$(x_{\theta}^-\otimes t^{(\lambda|\theta)+s^\prime+1})\,v=0\quad \textup{and}\quad (\theta^\vee\otimes t^{s\prime})\,v=0, \qquad \textup{for every}\,\, s^\prime\geq1.$$

We now prove part \eqref{lemmapart1b}. 
We observe  that $(x_{\theta}^-\otimes t^{(\lambda|\theta)+1})^{2k-\ell+1}\,v$ is an element of weight $\ell\lambda+(\ell-k-1)\theta$. 
Considering the $\mathfrak{sl}_2$  copy spanned by  $x_{\alpha}^+\otimes 1, x_{\alpha}^-\otimes 1$, and $\alpha^\vee\otimes 1$, the proof of 
part~\eqref{lemmapart1b}
follows by standard $\mathfrak{sl}_2$ calculations,  using part \eqref{lemmapart1a}.  The proof of part \eqref{lemmapart2}  is similar and we omit the details.
\end{proof}
The next proposition gives the existence of ${\phi^{-}_1}$ and ${\phi^{-}_2}$.
\begin{proposition}\label{existm}
Let $k, \ell, m\in\mathbb{N}$ and $\lambda\in L^+$. Denote $m^\prime=\ell+m-2k-1$ and $k^\prime=\ell-k-1$.
\be
\item\label{lemmapp11}
If $\ell\leq 2k$ and $k\leq m<\ell$, then the map $ \phi_2^-: \tau_{(2k-\ell+1)((\lambda|\theta)+1)}{\mathbf{V}}^{\ell,\,\ell\lambar}_{m^\prime,\,k^\prime\theta} \rightarrow \textup{ker}\, {\phi^{+}_m}$
which takes $v^{\ell,\,\ell\lambar}_{m^\prime,\,k^\prime\theta}\mapsto \big(x_\theta^{-}\otimes t^{(\lambda|\theta)+1}\big)^{2k-\ell+1}\,v^{\ell,\,\ell\lambar}_{m,\,k\theta}$ is a 
surjective morphism of $\lieg[t]$-modules.
\item\label{lemmapp22}
If $\ell\geq 2k$, then the map $\phi_1^-: \tau_{(\lambda|\theta)+1} {\mathbf{V}}^{\ell,\,\ell\lambar}_{k-1,\,(k-1)\theta} \rightarrow \textup{ker}\, {\phi^{+}_k} \,\,\text{which takes}\,\, v^{\ell,\,\ell\lambar}_{k-1,\,(k-1)\theta}\mapsto\\ \big(x^-_{\theta}\otimes t^{(\lambda|\theta)+1}\big)\,v^{\ell,\,\ell\lambar}_{k,\,k\theta}$ is a 
surjective morphism of $\lieg[t]$-modules.
\ee
\end{proposition}
\begin{proof}
Set $v=v^{\ell,\,\ell\lambar}_{m,\,k\theta}$.
To prove part \eqref{lemmapp11}, we need to show that $\big(x_\theta^{-}\otimes t^{(\lambda|\theta)+1}\big)^{2k-\ell+1}\,v$ satisfies the defining relations 
of the module ${\mathbf{V}}^{\ell,\,\ell\lambar}_{m^\prime,\,k^\prime\theta}$. Since $m^\prime\leq2k^\prime$,
using Lemma \ref{l1} \eqref{lemmapart1} and the relations \eqref{mr:dfr3}, we only need to show the following:
\begin{gather}
\big(x^{-}_{\alpha}\otimes t^{(\lambda|\alpha)}\big)^{\langle (k^\prime-i^\prime)\theta,\, \alpha^{\vee} \rangle +1} \big(x^{-}_{\theta}\otimes t^{(\lambda|\theta)+1}\big)^{i^\prime+2k-\ell+1}\,v=0, \quad \forall \,\,\alpha\in R^+,\,\,\, 0\leq i^\prime\leq k^\prime,\label{redrelinlem}
\\
\big(x_{\theta}^-\otimes t^{(\lambda|\theta)+1}\big)^{(2k^\prime-m^\prime+1)+(2k-\ell+1)}\,v=0.\label{mr:dfr5inlem}
\end{gather}
Since $k^\prime-i^\prime=k-(i^\prime+2k-\ell+1)$, the relations \eqref{redrelinlem} follow from the relations \eqref{redrel}. The relation \eqref{mr:dfr5inlem} is same as
the relation \eqref{mr:dfr5}. Hence part \eqref{lemmapp11}. The proof of part \eqref{lemmapp22}  is similar and we omit the details.
%The proof follows by using Lemma \ref{l1} and the relations \eqref{mr:dfr3}-\eqref{mr:dfr5}.
\end{proof}
\subsection{}
Using \eqref{lgeq2k}, Theorem \ref{CSVWT}, and Proposition \ref{Prop}, the existence of surjective map  ${\phi^{+}_m}$ and the maps ${\phi^{-}_1}$ and ${\phi^{-}_2}$ give the 
following inequalities:
\beq\label{diml1}
\dim \mathbf{V}^{\ell,\,\ell\lambar}_{k,\,k\theta} \leq \dim \mathbf{V}^{\ell,\,\ell\lambar}_{k-1,\,(k-1)\theta}+\dim V(k\theta) \prod_{j=1}^{p}{\dim D(\ell,\, \ell\lambda_j)},\, \,\, \forall \,\, \ell\geq 2k,
\eeq
\beq\label{diml2}
\dim \mathbf{V}^{\ell,\,\ell\lambar}_{m,\,k\theta} \leq \dim \mathbf{V}^{\ell,\,\ell\lambar}_{\ell+m-2k-1,\,(\ell-k-1)\theta}+\dim D(\ell,\,k\theta) \prod_{j=1}^{p} {\dim D(\ell,\, \ell\lambda_j)},\,\, \, \forall \,\, k\leq m<\ell\leq 2k,
\eeq
for every $\lambda_1,\ldots,\lambda_p\in {L^+}$ with $\lambda_1+\cdots+\lambda_p=\lambda$.

The following proposition is useful in getting the reverse inequalities.
\begin{proposition}\label{Vsurjfu}
Let $k, \ell, m\in\mathbb{N}$ be such that $\ell\geq m\geq k$. Let $\lambda_1,\ldots,\lambda_p\in {L^+}$ and denote $\lambda=\lambda_1+\cdots+\lambda_p$. 
Then the assignment $v^{\ell,\,\ell\lambar}_{m,\,k\theta} \mapsto w_{\ell,\,\ell\lambda_1}*\cdots*w_{\ell,\,\ell\lambda_p}*w_{m,\,k\theta}$ defines 
a surjective  morphism of $\lieg[t]$-modules
\[\mathbf{V}^{\ell,\,\ell\lambar}_{m,\,k\theta} \twoheadrightarrow 
D(\ell,\,\ell\lambda_1)*\cdots*D(\ell,\,\ell\lambda_p)*D(m,\,k\theta).\]
In particular,
$$
\dim \mathbf{V}^{\ell,\,\ell\lambar}_{m,\,k\theta} \geq  \dim D(m,\,k\theta)
\prod_{j=1}^{p} {\dim D(\ell,\, \ell\lambda_j)}.
$$
\end{proposition}
\begin{proof}
The proof  follows from  Lemma \ref{l2}. 
\end{proof}
\subsection{Proof of Theorem \ref{MTone}}
Let $\lambda_1,\ldots,\lambda_p\in {L^+}$ be such that $\lambda_1+\cdots+\lambda_p=\lambda$. 
\begin{enumerate} \item
To prove part \eqref{pone},  we prove that
\begin{equation}\label{mtpfp4eq1111}
{\mathbf{V}}^{\ell,\,\ell\lambda}_{k,\,k\theta}\cong_{\mathfrak{g}[t]}
D(\ell,\,\ell\lambda_1)*\cdots*D(\ell,\,\ell\lambda_p)*D(k,\,k\theta),\quad \textup{for}\,\,\ell\geq 2k.
\end{equation}
 We proceed by induction on $k$. In the case $k=1$, by using \eqref{obs1} and Theorem \ref{CSVWT}, we have
 $\dim \mathbf{V}^{\ell,\,\ell\lambar}_{0,\,0\theta}=\prod_{j=1}^{p} {\dim D(\ell,\, \ell\lambda_j)}.$
 Now  \eqref{diml1} gives,
 \beq\label{pfmtp1eq1}
 \dim \mathbf{V}^{\ell,\,\ell\lambar}_{1,\,1\theta} \leq (\dim V(\theta)+1) \prod_{j=1}^{p}{\dim D(\ell,\, \ell\lambda_j)} .\eeq
Using Proposition \ref{DlkTP}, we have 
\beq\label{pfmtp2eq2}\dim D(1,\,\theta)=\dim V(\theta)+1.\eeq
We get \eqref{mtpfp4eq1111}  in this case from Proposition \ref{Vsurjfu}, by using \eqref{pfmtp1eq1} and \eqref{pfmtp2eq2}.

Now suppose $k\geq 2$. By the induction hypothesis, we have
\[\dim \mathbf{V}^{\ell,\,\ell\lambar}_{k-1,\,(k-1)\theta}=\dim D(k-1,\,(k-1)\theta) \prod_{j=1}^{p} {\dim D(\ell,\, \ell\lambda_j)}.\]
Substituting this into \eqref{diml1}, we get
\[\dim \mathbf{V}^{\ell,\,\ell\lambar}_{k,\,k\theta} \leq \big(\dim V(k\theta)+\dim D(k-1,\,(k-1)\theta)\big) \prod_{j=1}^{p}{\dim D(\ell,\, \ell\lambda_j)} .\]
%Proposition \ref{DlkTP} gives,
%\[\dim \mathbf{V}^{\ell,\,\ell\lambar}_{k,\,k\theta} \leq \dim D(k,\,k\theta) \prod_{j=1}^{p} {\dim D(\ell,\, \ell\lambda_j)}.\]
Now the proof of \eqref{mtpfp4eq1111} follows from Proposition \ref{Vsurjfu}, by using Proposition \ref{DlkTP}.
\item
The existence of ${\psi}^{-}$ such that
\[{\psi}^{-}\big({v}^{\ell,\,\ell\lambar}_{m-k-1,\,(m-k-1)\theta}\big)= \big(x^{-}_{\theta}\otimes t^{(\lambda|\theta)+1}\big)^{2k-m+1}\,{v}^{\ell,\,\ell\lambar}_{k,\,k\theta},\]
and ${\psi}^{-}$ is injective  follow from  repeated applications of the injective morphism ${\phi^{-}_1}$ from  part \eqref{pone} for various values of $k$.
It is easy to see that the assignment ${v}^{\ell,\,\ell\lambar}_{k,\,k\theta} \mapsto {v}^{\ell,\,\ell\lambar}_{m,\,k\theta}$ gives the existence  of ${\psi}^+$. Clearly ${\psi}^+$ is 
surjective and 
\[\textup{ker}\,{\psi}^+=\mathbf{U}(\lieg[t])\,\big(x^{-}_{\theta}\otimes t^{(\lambda|\theta)+1}\big)^{2k-m+1}\,{v}^{\ell,\,\ell\lambar}_{k,\,k\theta}.\] 
Now we  observe that $\textup{img}\,{\psi}^{-}=\textup{ker}\,{\psi}^+$.
This completes the proof of part \eqref{ptwo}, and it gives,
\[\dim\,\mathbf{V}^{\ell,\,\ell\lambar}_{m,\,k\theta}=\dim\,\mathbf{V}^{\ell,\,\ell\lambar}_{k,\,k\theta}-\dim\,\mathbf{V}^{\ell,\,\ell\lambar}_{m-k-1,\,(m-k-1)\theta}.\]
Since $\ell\geq 2k> 2(m-k-1)$, using \eqref{mtpfp4eq1111}, the right hand side of
the last equation becomes
\[\big(\dim D(k,\,k\theta)-\dim D({m-k-1,\,(m-k-1)\theta})\big) \prod_{j=1}^{p} {\dim D(\ell,\, \ell\lambda_j)}.\]
Now using Propositions \ref{DlkTP} and \ref{Vsurjfu}, we get
\begin{equation}\label{mtpfp4eq1112}
{\mathbf{V}}^{\ell,\,\ell\lambda}_{m,\,k\theta}\cong_{\mathfrak{g}[t]}
D(\ell,\,\ell\lambda_1)*\cdots*D(\ell,\,\ell\lambda_p)*D(m,\,k\theta),\quad\textup{for}\,\,k<m\leq 2k\leq\ell.
\end{equation}
\item 
Let $k^\prime=\ell-k-1$ and $m^\prime=\ell+m-2k-1$.
Since $\ell\geq 2k^\prime$ and $k^\prime\leq m^\prime \leq 2k^\prime$,
using \eqref{mtpfp4eq1111} and \eqref{mtpfp4eq1112}, we get
\[\dim \mathbf{V}^{\ell,\,\ell\lambar}_{m^\prime,\,k^\prime\theta}= \dim D(m^\prime,\,k^\prime\theta) \prod_{j=1}^{p} {\dim D(\ell,\, \ell\lambda_j)}.\]
Substituting this  into \eqref{diml2}, we get 
\[\dim \mathbf{V}^{\ell,\,\ell\lambar}_{m,\,k\theta} \leq \big(\dim D(m^\prime,\,k^\prime\theta)+\dim D(\ell,\,k\theta)\big) \prod_{j=1}^{p} {\dim D(\ell,\, \ell\lambda_j)} .\]
Now Proposition \ref{DlkTP} gives,
\[\dim \mathbf{V}^{\ell,\,\ell\lambar}_{m,\,k\theta} \leq \dim D(m,\,k\theta)\prod_{j=1}^{p} {\dim D(\ell,\, \ell\lambda_j)}.\]
Using Proposition \ref{Vsurjfu}, we get part \eqref{pthree}, and
\begin{equation}\label{mtpfp4eq1113}
{\mathbf{V}}^{\ell,\,\ell\lambda}_{m,\,k\theta}\cong_{\mathfrak{g}[t]}
D(\ell,\,\ell\lambda_1)*\cdots*D(\ell,\,\ell\lambda_p)*D(m,\,k\theta), \quad\textup{for}\,\,k\leq m<\ell\leq 2k.
\end{equation}
\end{enumerate}
\qed
\subsection{}
%In this subsection, we discuss about the generalized Demazure modules.
Although the following proposition seems to be well known, we give a proof here for the sake of completeness.
\begin{proposition}\label{gdmvsgt}
For $1\leq j\leq p$, let $ \Lambda^j\in\widehat{P}^+,  \xi_j\in\afw \Lambda^j$  such that $\langle \xi_j,\, \alpha_i^{\vee}\rangle \leq 0, \,\forall\, i\in I$, and $v_{w_0\xi_j}\in V(\Lambda^j)_{w_0 \xi_j}$. Then
\[D\big(\xi_1,\ldots,\xi_p\big)=\mathbf{U}(\lieg[t])\,(v_{w_0 \xi_1}\otimes\cdots\otimes v_{w_0\xi_p}).\]
\end{proposition}
\begin{proof}
We first prove that $v_{w_0 \xi_1}\otimes\cdots\otimes v_{w_0 \xi_p}\in D\big(\xi_1,\ldots,\xi_p\big)$. Let 
$v_{\xi_j}\in V(\Lambda^j)_{\xi_j}$, for $1\leq j\leq p$, and set $v=v_{\xi_1}\otimes\cdots\otimes v_{\xi_p}.$ Since $\mathfrak{n}^-\,v=0$, the $\lieg$-submodule 
$\mathbf{U}(\lieg)\,v\subseteq D\big(\xi_1,\ldots,\xi_p\big)$ is isomorphic to the irreducible highest weight $\lieg$-module $V\big(w_0(\xi_1+\cdots+\xi_p)|_{\csa}\big)$. 
Hence there exists a non-zero element $v^{\prime}\in\mathbf{U}(\mathfrak{n}^+)\,v$ in 
$D\big(\xi_1,\ldots,\xi_p\big)$, whose $\csa$-weight is equal to $w_0(\xi_1+\cdots+\xi_p)|_{\csa}$. Suppose that $\xi_1+\cdots+\xi_p=a\Lambda_0+\beta+b\delta,$ for some $a, b\in\complex$ and  
$\beta\in{\csa}^{*}$, then the weight of $v^{\prime}$ is equal to $ a\Lambda_0+w_0\beta+b\delta=w_{0}(\xi_1+\cdots+\xi_p).$
Hence $v^{\prime}$ is a non-zero constant multiple of $v_{w_0\xi_1}\otimes\cdots\otimes v_{w_0 \xi_p}$, since
\[\big(D(\xi_1)\otimes\cdots\otimes D(\xi_p)\big)_{w_0(\xi_1+\cdots+\xi_p)}=\complex\,(v_{w_0 \xi_1}\otimes\cdots\otimes v_{w_0 \xi_p}).\]
Now we have
\[\mathbf{U}(\lieg[t])\,(v_{w_0 \xi_1}\otimes\cdots\otimes v_{w_0\xi_p})\subseteq D\big(\xi_1,\ldots,\xi_p\big).\]
The reverse containment also follows in the similar way and we omit the details.
\end{proof}
%The following proposition  is useful.
\begin{proposition}\label{p}
Let $k, \ell, m\in\mathbb{N}$ be such that $\ell\geq m\geq k$ and $\lambda\in L^+$. 
%Let $\lambar=(\lambda_1,\ldots,\lambda_p)\in ({L^+})^{p}$, for $p\in\mathbb{N}$, and denote that $\lambda=\lambda_1+\cdots+\lambda_p$. 
Then there exists a morphism of $\lieg[t]$-modules
\[\varphi: \mathbf{V}^{\ell,\,\ell\lambar}_{m,\,k\theta}\rightarrow D(\ell-m,\,(\ell-m)\lambda)\otimes D(m,\,m\lambda+k\theta),\]
such that
\[\varphi \big({v}^{\ell,\,\ell\lambar}_{m,\,k\theta})= w_{\ell-m,\,(\ell-m)\lambda}\otimes w_{m,\,m\lambda+k\theta} \qquad\text{and} \qquad\textup{img}\,\varphi=\mathbf{U}(\lieg[t])\,(w_{\ell-m,\,(\ell-m)\lambda}\otimes w_{m,\,m\lambda+k\theta}).\]
\end{proposition}
\begin{proof}
The proof follows by using  Proposition \ref{Dllamkt} and the relations 
$$(x_{\alpha}^-\otimes t^{(\lambda|\alpha)})\,w_{\ell-m,\,(\ell-m)\lambda}=0, \quad\forall\,\,\alpha\in R^+,$$
which hold in the module $D(\ell-m,\,(\ell-m)\lambda)$.
\end{proof}
%The following proposition helps in proving our main result (Theorem \ref{MT1}).
\begin{proposition}\label{VsurjGD}
With hypothesis and notation as in Proposition \ref{p}, there exist  surjective morphisms of $\lieg[t]$-modules,
\begin{gather*}
\mathbf{V}^{\ell,\,\ell\lambar}_{m,\,k\theta}\twoheadrightarrow \begin{cases} D\big(t_{w_0\lambda}(\ell-m)\Lambda_0,\, t_{w_0 (\lambda+\theta)}  
(m\Lambda_0+(m-k)\theta)\big),  &   m\leq 2k,
\\
D\big(t_{w_0\lambda}(\ell-m)\Lambda_0,\, t_{w_0\lambda}w_0(m\Lambda_0+k\theta)\big), &  m\geq 2k.\end{cases}
\end{gather*}
\end{proposition}
\begin{proof}
We observe from Proposition \ref{DvsV} that
\[
D(\ell-m,\,(\ell-m)\lambda)\cong_{\lieg[t]}
D(t_{w_0\lambda}(\ell-m)\Lambda_0)\]
and
\[D(m,\,m\lambda+k\theta)\cong_{\lieg[t]}\begin{cases} D(t_{w_0(\lambda+\theta)}(m\Lambda_0+(m-k)\theta)), &   m\leq 2k,\\
                                          D(t_{w_0\lambda}w_0(m\Lambda_0+k\theta)), &  m\geq 2k.
                                         \end{cases}
\]
Now the proof follows from Proposition \ref{p}, by using Proposition \ref{gdmvsgt}.
\end{proof}
\subsection{}
In this subsection, we recollect some facts about the Demazure operators and the character of generalized Demazure modules, which are useful in proving Theorem \ref{MTtwo}.

For $0\leq i \leq n$, the Demazure operator $\mathcal{D}_i$ is a linear operator on $\mathbb{Z}[\widehat{P}]$, and is defined by
\[\mathcal{D}_i(e^{\Lambda})= \frac{e^\Lambda-e^{r_i(\Lambda)-\alpha_i}}{1-e^{-\alpha_i}}.\]
For $w\in\widehat{W}$ and a reduced expression $w=r_{i_1}\cdots r_{i_k}$, the Demazure operator $\mathcal{D}_w$ is defined as 
$\mathcal{D}_w=\mathcal{D}_{i_1}\cdots \mathcal{D}_{i_k}$, and is independent of the choice of reduced expression of $w$ (\cite[Corollary~8.2.10]{kum2}).  
For $w\in\widehat{W}$ and $\sigma\in\Sigma$, set $\mathcal{D}_{w\sigma}(e^\Lambda)=\mathcal{D}_{w}(e^{\sigma\Lambda}).$ Since  $\mathcal{D}_i(e^\delta)=e^\delta,$ the operator $\mathcal{D}_w$ descends to $\mathbb{Z}[\widehat{P}]/I_{\delta},$ for all $w\in\widetilde{W}.$

The following theorem gives the character of generalized Demazure modules in terms of Demazure operators. It is a combination of the results \cite[Proposition~2.7 and Corollary 2.8]{Naoi3}. The key ingredient in the proof is a result from \cite{LLM}.
\begin{theorem}\label{naoi}\cite{LLM, Naoi3}
Let $\Lambda^1,\ldots,\Lambda^p$ be a sequence of elements of $\widehat{P}^+$. Let $w_1,\ldots,w_p$ be a sequence of elements of $\widetilde{W}$, and denote $w_{[1,j]}=w_1\cdots w_j,\, \forall\, 1\leq j\leq p.$ If $\ell(w_{[1,p]})=\sum_{j=1}^{p} \ell(w_j)$, then we have 
\begin{equation*}
\begin{split}
&\textup{ch}_{\widehat{\mathfrak{h}}}\,D\left(w_{[1,1]}\Lambda^1, w_{[1,2]}\Lambda^2, \ldots, w_{[1,p-1]}\Lambda^{p-1}, w_{[1,p]}\Lambda^p\right)\\
&=\mathcal{D}_{w_1}\left(e^{\Lambda^1}\,\mathcal{D}_{w_2}\left(e^{\Lambda^2} \cdots \mathcal{D}_{w_{p-1}}\left(e^{\Lambda^{p-1}}\,\mathcal{D}_{w_p}(e^{\Lambda^p})\right)\cdots\right)\right).
\end{split}
\end{equation*}
\end{theorem}
The following theorem may be found in \cite[Theorem 3.4]{kum1} and \cite{M} (see also \cite[Theorem 8.2.9]{kum2} and \cite[\S {4.5}]{CSVW}).
\begin{theorem}\cite{kum1, M}\label{CVK}
Given $(\ell,\lambda)\in\mathbb{N}\times P^+$, let $w\in\widehat{W}, \sigma\in\Sigma$ and $\Lambda\in \widehat{P}^+$ such that 
\[w\sigma\Lambda\equiv w_0\lambda+\ell\Lambda_0 \,\,\,\textup{mod} \,\mathbb{Z}\delta.\] Then
\[\mathcal{D}_{w\sigma}(e^\Lambda)\equiv e^{\ell\Lambda_0} \textup{ch}_{\csa}\,D(\ell,\lambda)\,\,\,\textup{mod}\,I_{\delta}.\]
\end{theorem}
The next lemma follows as \cite[Lemma 7]{FoL1} (see also \cite[\S {4.4}]{CSVW}). 
\begin{lemma}\cite{FoL1}\label{CV}
Let $V$ be a finite-dimensional $\lieg$-module. Let $(\ell, \mu)\in\mathbb{N}\times L^+$. Then for $A\in\mathbb{Z}$, we have
\[\mathcal{D}_{t_{w_0\mu}}(e^{\ell\Lambda_0+A\delta}\,\textup{ch}_{\csa}\,V)= \mathcal{D}_{t_{w_0\mu}}(e^{\ell\Lambda_0+A\delta})\, \textup{ch}_{\csa}\,V.\]
\end{lemma}
\subsection{}
We now prove that the fusion product of  Demazure modules of different level as a $\lieg$-module is isomorphic to a generalized Demazure module. We further conjecture that they are in fact isomorphic as $\lieg[t]$-modules.
%The following proposition helps in proving our main result.
\begin{proposition}\label{hiso}
Let $(\ell,\lambda)\in \mathbb{N}\times P^+$. Suppose that there exist $\nu\in P^+, w\in W$, and $\Lambda\in \widehat{P}^+$  such that 
\[t_{-\nu}w\Lambda\equiv w_0\lambda+\ell\Lambda_0\,\,\,\textup{mod}\,\, \mathbb{Z}\delta.\]
Let $(\ell_i,\lambda_i)\in  \mathbb{N}\times L^+$, for $1\leq i\leq p$, be such that $\ell_1\geq \cdots\geq \ell_p\geq \ell$.  
Then we have the following isomorphism as $\lieg$-modules:
\begin{gather}
\begin{split}
&D(\ell_1,\,\ell_1\lambda_1)*\cdots * D(\ell_p,\,\ell_p\lambda_p) * D(\ell,\,\lambda)
%\nonumber
\\
&\cong D\big( t_{w_0\lambda_1}(\ell_1-\ell_2)\Lambda_0,\ldots,t_{w_0(\lambda_1+\cdots+\lambda_p)}(\ell_p-\ell)\Lambda_0,\,t_{w_0(\lambda_1+\cdots+\lambda_p)}\,t_{-\nu}w\Lambda\big).\label{gconj}
\end{split}\end{gather}
\end{proposition}
\begin{proof}
Since the finite-dimensional $\lieg$-modules are determined by their $\csa$-characters, it suffices to show that
\begin{equation}
\begin{split}
&\textup{ch}_{\csa}\,D\big( t_{w_0\lambda_1}(\ell_1-\ell_2)\Lambda_0,\ldots,\,t_{w_0(\lambda_1+\cdots+\lambda_p)}(\ell_p-\ell)\Lambda_0,\,t_{w_0(\lambda_1+\cdots+\lambda_p)}\,t_{-\nu}w\Lambda\big)\\
&=\textup{ch}_{\csa}\,D(\ell_1,\,\ell_1\lambda_1)\cdots\textup{ch}_{\csa}\,D(\ell_p,\,\ell_p\lambda_p)\,\textup{ch}_{\csa}D(\ell,\,\lambda).\label{hchar}
\end{split}\end{equation}
From  Lemma \ref{lengthl} and Theorem \ref{naoi},  we have
\begin{equation*}
\begin{split}
&\textup{ch}_{\widehat{\csa}}\,D\big( t_{w_0\lambda_1}(\ell_1-\ell_2)\Lambda_0,\ldots,\,t_{w_0(\lambda_1+\cdots+\lambda_p)}(\ell_p-\ell)\Lambda_0,\,t_{w_0(\lambda_1+\cdots+\lambda_p)}\,t_{-\nu}w\Lambda\big)
\\
&=\mathcal{D}_{t_{w_0\lambda_1}}\big(e^{(\ell_1-\ell_2)\Lambda_0}\,\mathcal{D}_{t_{w_0\lambda_2}}\big(e^{(\ell_2-\ell_3)\Lambda_0}\cdots \mathcal{D}_{t_{w_0\lambda_p}}\big(e^{(\ell_p-\ell)\Lambda_0}\, \mathcal{D}_{t_{-\nu}w}\big(e^{\Lambda}\big)\big)\cdots\big)\big).
\end{split}
\end{equation*}
Now repeated applications of Theorem \ref{CVK} and Lemma \ref{CV} give,
\begin{gather*}
\begin{split}
&\textup{ch}_{\widehat{\csa}}\,D\big( t_{w_0\lambda_1}(\ell_1-\ell_2)\Lambda_0,\ldots,\,t_{w_0(\lambda_1+\cdots+\lambda_p)}(\ell_p-\ell)\Lambda_0,\,t_{w_0(\lambda_1+\cdots+\lambda_p)}\,t_{-\nu}w\Lambda\big)\\
&\equiv e^{\ell_1\Lambda_0}\,\textup{ch}_{\csa}\,D(\ell_1,\,\ell_1\lambda_1)\cdots\textup{ch}_{\csa}\,D(\ell_p,\,\ell_p\lambda_p)\,\textup{ch}_{\csa}D(\ell,\,\lambda) \quad\,\,\textup{mod}\,\,I_{\delta}.
\end{split}
\end{gather*}
Letting $e^{\Lambda_0}\mapsto 1$ and $e^\delta\mapsto 1,$ we get \eqref{hchar}. Hence the proposition.
\end{proof}
We conjecture below that the isomorphism \eqref{gconj} also holds as $\lieg[t]$-modules.
\begin{conjecture}\label{conj}
Under the hypothesis of Proposition \ref{hiso}, we have the isomorphism \eqref{gconj} of $\lieg[t]$-modules.
\end{conjecture} 
\begin{remark}
This conjecture is proved when $\ell_1=\cdots=\ell_p=\ell$ in \cite{CSVW}.  
In \cite{Naoi2}, for $\lieg$ simply laced, it is proved  when $\lambda_1,\ldots,\lambda_p$, and $\frac{1}{\ell}\lambda$ are fundamental weights. 
In this paper, we prove  a special case of it and also give the defining relations (see Theorem \ref{MTtwo}).
\end{remark}
\subsection{Proof of Theorem \ref{MTtwo}}
We first prove that for $\ell\geq m\geq k$,
\begin{equation}\label{mtpfp4eq1}
D(\ell,\,\ell\lambda_1)*\cdots*D(\ell,\,\ell\lambda_p)*D(m,\,k\theta)\cong_{\mathfrak{g}[t]}{\mathbf{V}}^{\ell,\,\ell\lambda}_{m,\,k\theta}.
\end{equation}
We consider two cases. First, suppose that $\ell\geq 2k$. Then, from \eqref{mtpfp4eq1111} and \eqref{mtpfp4eq1112}  we obtain \eqref{mtpfp4eq1} when $k\leq m\leq 2k$. If $m>2k$, then since 
$D(m,\,k\theta)\cong_{\lieg[t]} D(2k,\,k\theta)$,
it follows from the case $m=2k$ and using \eqref{obs0}. In the second case, we have $\ell\leq 2k$. Then, from \eqref{mtpfp4eq1113} we obtain \eqref{mtpfp4eq1} when $\ell\neq m$. If $\ell=m$, then
it follows by using Theorem \ref{CSVWT} and Proposition \ref{Prop}.

Next we prove that for $\ell\geq m\geq k$,
\beq\label{mtpfp4eq11}
{\mathbf{V}}^{\ell,\,\ell\lambar}_{m,\,k\theta}
\cong_{\lieg[t]}\begin{cases} D\big(t_{w_0\lambda} (\ell-m)\Lambda_0,\, t_{w_0(\lambda+\theta)} (m\Lambda_0+(m-k)\theta)\big), &   m\leq 2k,\\
D\big(t_{w_0\lambda}(\ell-m)\Lambda_0,\, t_{w_0\lambda}w_0 (m\Lambda_0+k\theta)\big), &  m\geq 2k.\end{cases}\eeq
Using \eqref{DisoV} and \eqref{lgeq2k}, Proposition \ref{hiso} gives,
\beq\label{mtpfp4eq2}
\begin{split}
&\dim D(\ell,\, \ell\lambda) \dim D(m,\,k\theta)\\&=
\begin{cases}
\dim D\big(t_{w_0\lambda} (\ell-m)\Lambda_0,\, t_{w_0(\lambda+\theta)} (m\Lambda_0+(m-k)\theta)\big), &   m\leq 2k,\\
\dim D\big(t_{w_0\lambda}(\ell-m)\Lambda_0,\, t_{w_0\lambda}w_0 (m\Lambda_0+k\theta)\big), &  m\geq 2k.\end{cases}\end{split}
\eeq
Since $\lambda=\lambda_1+\cdots+\lambda_p$, we have from Theorem \ref{CSVWT} that
\beq\label{mtpfp4eq3}
\dim D(\ell,\, \ell\lambda)=\prod_{j=1}^{p} {\dim D(\ell,\, \ell\lambda_j)}.
\eeq
The proof of \eqref{mtpfp4eq11} follows from Proposition \ref{VsurjGD}, by using \eqref{mtpfp4eq1}, \eqref{mtpfp4eq2}, and \eqref{mtpfp4eq3}.
\qed
%We now prove  part (2) of Theorem \ref{MT1}. 
%We now prove part (3) of Theorem \ref{MT1}. 
%We now prove Part (4) of Theorem \ref{MT1}. 
\section{The connection with Chari-Venkatesh modules}\label{connectiontoCV}
In this section, we prove that the defining relations of $\mathbf{V}^{\ell,\,\ell\lambar}_{m,\,k\theta}$ can be simplified. 
This allows us to make connection with the modules introduced by 
Chari and Venkatesh in \cite{CV}.
\subsection{}
Let us begin with recalling the notations and definition given in \cite{CV}. 
%In this subsection we recall the modules introduced in \cite{CV}. 
%Following \cite{CV}, we first introduce some notation. 
For $r, s\in\mathbb{Z}_{\geq 0}$, let
\[\mathbf{S}(r,s)=\{(b_p)_{p\geq0} : b_p\in\mathbb{Z}_{\geq0},\, \sum_{p\geq 0}b_p =r,\, \sum_{p\geq0}pb_p =s\}.\]
For $\alpha\in R^{+}$ and $r,s\in\mathbb{Z}_{\geq0}$, define an element $\mathbf{x}^{-}_{\alpha}(r,s)\in\mathbf{U}(\lieg[t])$ by
\begin{equation}\label{prelim}
\mathbf{x}^{-}_{\alpha}(r,s)=\sum_{(b_p)\in\mathbf{S}(r,s)} (x^{-}_{\alpha} \otimes 1)^{(b_0)}(x^{-}_{\alpha} \otimes t)^{(b_1)}\cdots(x^{-}_{\alpha} \otimes t^s)^{(b_s)},
\end{equation}
where  $X^{(b)}$ denotes the divided power $X^b/b!$.
For $\kay\in\mathbb{Z}_{\geq0}$, let $ {\mathbf{S}}(r,s)_\kay$ (resp. $_{\kay}{\mathbf{S}}(r,s)$) be the subset of $\mathbf{S}(r,s)$ consisting of elements $(b_p)_{p\geq0},$ 
satisfying $b_p=0$ for $p\geq \kay$ (resp. $b_p=0$ for $p<\kay$).
For $\alpha\in R^{+}$ and $r,s, \kay\in\mathbb{Z}_{\geq0}$, define elements $\mathbf{x}^{-}_{\alpha}(r,s)_\kay$ and  $_\kay{\mathbf{x}^{-}_{\alpha}}(r,s)$ of $\mathbf{U}(\lieg[t])$
by
\begin{gather}
{\mathbf{x}^{-}_{\alpha}}(r,s)_\kay=\sum_{(b_p)\in{\mathbf{S}}(r,s)_\kay}(x^{-}_{\alpha} \otimes 1)^{(b_0)}(x^{-}_{\alpha} \otimes t^{1})^{(b_{1})}\cdots(x^{-}_{\alpha} \otimes t^{\kay-1})^{(b_{\kay-1})},\label{prelimkr}
\\
_\kay{\mathbf{x}^{-}_{\alpha}}(r,s)=\sum_{(b_p)\in_\kay{\mathbf{S}}(r,s)} (x^{-}_{\alpha} \otimes t^\kay)^{(b_\kay)}(x^{-}_{\alpha} \otimes t^{\kay+1})^{(b_{\kay+1})}\cdots(x^{-}_{\alpha} \otimes t^s)^{(b_s)}.\label{prelimk}
\end{gather}

Given $\mu\in P^+$ and a $|R^+|$-tuple $\bm{\xi}=(\xi(\alpha))_{\alpha\in R^+}$ of partitions such that $|\xi(\alpha)|=\langle\mu,\,\alpha^{\vee}\rangle,\,\,\forall 
\,\,\alpha\in R^{+}$. The module $V(\bm{\xi})$ is the cyclic $\lieg[t]$-module generated by $v_{\bm{\xi}}$ with defining relations:
\begin{gather}
(x^{+}_{i}\otimes t^s)\,v_{\bm{\xi}}=0,\quad (\hi\otimes t^s)\,v_{\bm{\xi}} =\delta_{s,0} \langle \mu, \, \hi \rangle\, v_{\bm{\xi}},  \quad \big(x^{-}_{i}\otimes 1\big)^{\langle \mu, \, \hi \rangle +1} 
\,v_{\bm{\xi}}=0,\quad\forall \,\, s\geq0, i\in I,\label{cvs:dfr1}
%\mathfrak{n}^+[t]\,v_{\bm{\xi}}=0,\qquad (h\otimes t^s)\,v_{\bm{\xi}}=\delta_{s,0} \langle \mu, \, h \rangle\, v_{\bm{\xi}}, 
%\quad \forall \,\, h\in\csa, s\in\mathbb{Z}_{\geq0},\label{cvs:dfr1}
%\\
%\big(x^{-}_{\alpha}\otimes 1\big)^{\langle \mu, \, \alpha^{\vee} \rangle +1} \,v_{\bm{\xi}}=0,\quad \forall \,\, \alpha\in R^+,\label{cvs:dfr2}
\\
\mathbf{x}^{-}_{\alpha}(r,s)\,v_{\bm{\xi}}=0,\quad\forall\,\, \alpha\in R^+,  s, r \in \mathbb{N} \,\, \textup{such that} \,\, s+r\geq 1+r\kay+\sum_{j\geq \kay+1} \xi(\alpha)_j \,\,
\textup{for some} \,\, \kay\in \mathbb{N}.\label{cvs:dfr3}
\end{gather}
It is proved in \cite{CV} that the relations \eqref{cvs:dfr3} may be replaced with the following:
\beq\label{cvs:dfr4}
_\kay{\mathbf{x}^{-}_{\alpha}}(r,s)\,v_{\bm{\xi}}=0,\quad\forall\,\, \alpha\in R^+,  s, r, \kay \in \mathbb{N} \,\, \textup{such that} \,\, s+r\geq 1+r\kay+\sum_{j\geq \kay+1} \xi(\alpha)_j.\eeq
\subsection{}
%In this subsection, we state our main result of this section and discuss applications.
For $k, \ell, m\in\mathbb{N}$ such that $\ell\geq m\geq k$ and $\lambda\in {L^+}$,
%Let $\lambar=(\lambda_1,\ldots,\lambda_p)\in (L^+)^p$ and denote $\lambda=\sum_{i=1}^p\lambda_i$. 
we define three $|R^{+}|$-tuple of partitions as follows:
\begin{align}
&\bm{\xi}(\ell,\,\ell\lambda):=\left({\xi}(\ell,\,\ell\lambda)(\alpha)\right)_{\alpha\in R^+},  \quad
\text{where}
\quad
{\xi}(\ell,\,\ell\lambda)(\alpha):=\left((d_\alpha\ell)^{(\lambda|\alpha)}\right),\nonumber
\\
&\bm{\xi}(m,\,k\theta):=\left({\xi}(m,\,k\theta)(\alpha)\right)_{\alpha\in R^+},  \quad
\text{where}
\quad
%\begin{align*}
\xi({m,\,k\theta})(\alpha):=
\begin{cases}
\emptyset, & (\theta|\alpha)=0,
\\
\left(d_\alpha k\right), & (\theta|\alpha)=1,
\\
\left(m,\,2k-m\right), & \alpha=\theta \,\,\, \textup{and}\,\,\, m\leq2k,
\\
\left(2k\right), & \alpha=\theta \,\,\, \textup{and}\,\,\, m\geq2k,                                                                                           
\end{cases}\nonumber
\\
&\bm{\xi}^{\ell,\,\ell\lambar}_{m,\,k\theta}:=\left(\xi^{\ell,\,\ell\lambar}_{m,\,k\theta}(\alpha)\right)_{\alpha\in R^+},  \quad
\text{where}
\quad
%\begin{align*}
\xi^{\ell,\,\ell\lambar}_{m,\,k\theta}(\alpha):=
\begin{cases}
\left((d_\alpha\ell)^{(\lambda|\alpha)}\right), &(\theta|\alpha)=0,
\\
\left((d_\alpha\ell)^{(\lambda|\alpha)},\,d_\alpha k\right), &(\theta|\alpha)=1,
\\
\left(\ell^{(\lambda|\theta)},\,m,\,2k-m\right), &\alpha=\theta \,\,\, \textup{and}\,\,\, m\leq2k,
\\
\left(\ell^{(\lambda|\theta)},\,2k\right), &\alpha=\theta \,\,\, \textup{and}\,\,\, m\geq2k.                                                                                             
\end{cases}\label{xilmk}
%\end{align*}
\end{align}
The following isomorphisms follow from \cite[Theorem 2]{CV}:
\beq\label{isocv1}
V(\bm{\xi}(\ell,\,\ell\lambda))\cong_{\lieg[t]} D(\ell,\,\ell\lambda)
\qquad\textup{and}\qquad
V(\bm{\xi}(m,\,k\theta))\cong_{\lieg[t]} D(m,\,k\theta).
\eeq

%The following theorem is the main result of this section.
We are now in position to state the main result of this section.
\begin{theorem}\label{MTCV}
Let $k, \ell, m\in\mathbb{N}$ be such that $\ell\geq m\geq k$ and $\lambda\in L^+$.
 Then there exists an isomorphism of $
\lieg[t]$-modules,
$$V(\bm{\xi}^{\ell,\,\ell\lambar}_{m,\,k\theta})\cong\mathbf{V}^{\ell,\,\ell\lambar}_{m,\,k\theta}.$$
\end{theorem}
\begin{remark}
 For $\mathfrak{g}=\mathfrak{sl}_2$, it is proved in \cite[\S6]{CV} that 
 the modules $V(\bm{\xi})$ 
 are fusion products of evaluation modules $\textup{ev}_0 \,V(r\varpi_1), r\in\mathbb{Z}_{\geq 0}$ and vice-versa. 
 Using this, Theorem \ref{MTCV} gives us  for $n\in\mathbb{Z}_{\geq0}$ that 
 $$\mathbf{V}^{\ell,\,\ell n\varpi_1}_{m,\,k\theta}\cong_{\mathfrak{sl}_2[t]}\begin{cases}
                                                                            \big(\textup{ev}_0\, V(\ell\varpi_1)\big)^{*n}*\textup{ev}_0 \,V(m\varpi_1)*\textup{ev}_0 \,V((2k-m)\varpi_1), &m\leq 2k,\\
                                                                            \big(\textup{ev}_0\, V(\ell\varpi_1)\big)^{*n}*\textup{ev}_0\, V(2k\varpi_1), &m\geq 2k.
                                                                          \end{cases}$$
                                                                          \end{remark}

%Theorem \ref{MTCV} is proved in \S\ref{ss:pf}.
The following corollary shows that certain types of generalized Demazure modules also belong to the family of modules defined in \cite{CV}.
\begin{corollary}\label{cor:all}
Let $k, \ell, m\in\mathbb{N}$ be such that $\ell\geq m\geq k$. Let  $\lambda_1,\ldots,\lambda_p$ be a sequence of elements of ${L^+}$, and denote $\lambda=\lambda_1+\cdots+\lambda_p$. 
Then 
%With hypothesis and notation as in Theorem \ref{MTCV}, 
we have the following isomorphisms of $
\lieg[t]$-modules:
\begin{gather*}
\begin{split}
&V(\bm{\xi}(\ell,\,\ell\lambda_1))*\cdots*V(\bm{\xi}(\ell,\,\ell\lambda_p))*V(\bm{\xi}(m,\,k\theta))\\
&\cong
D(\ell,\,\ell\lambda_1)*\cdots*D(\ell,\,\ell\lambda_p)*D(m,\,k\theta)\\
&\cong{\mathbf{V}}^{\ell,\,\ell\lambar}_{m,\,k\theta}
\cong V(\bm{\xi}^{\ell,\,\ell\lambar}_{m,\,k\theta})
\cong \begin{cases} D\big(t_{w_0\lambda} (\ell-m)\Lambda_0,\, t_{w_0(\lambda+\theta)} (m\Lambda_0+(m-k)\theta)\big), &  m\leq 2k,\\
D\big(t_{w_0\lambda}(\ell-m)\Lambda_0,\, t_{w_0\lambda}w_0 (m\Lambda_0+k\theta)\big), &  m\geq 2k.\end{cases}
\end{split}
\end{gather*}
%where $\lambda=\lambda_1+\cdots+\lambda_p$.
\end{corollary}
\begin{proof}
The proof is immediate from   Theorem \ref{MTtwo}, \eqref{isocv1}, and Theorem \ref{MTCV}. 
\end{proof}
\subsection{{{Proof of Theorem \ref{MTCV}}.}}
The proof  is immediate from the following two propositions.
\begin{proposition}\label{prop1}
With hypothesis and notation as in Theorem \ref{MTCV}, 
the assignment $v_{\bm{\xi}^{\ell,\,\ell\lambar}_{m,\,k\theta}}\mapsto v^{\ell,\,\ell\lambar}_{m,\,k\theta}$ gives a surjective
morphism from 
 $V(\bm{\xi}^{\ell,\,\ell\lambar}_{m,\,k\theta})$ onto $\mathbf{V}^{\ell,\,\ell\lambar}_{m,\,k\theta}.$
\end{proposition}
\begin{proposition}\label{prop2}
With hypothesis and notation as in Theorem \ref{MTCV}, the assignment $v^{\ell,\,\ell\lambar}_{m,\,k\theta}\mapsto v_{\bm{\xi}^{\ell,\,\ell\lambar}_{m,\,k\theta}}$ gives a surjective
morphism from 
 $\mathbf{V}^{\ell,\,\ell\lambar}_{m,\,k\theta}$ onto $V(\bm{\xi}^{\ell,\,\ell\lambar}_{m,\,k\theta})$.
\end{proposition}
\subsubsection{{{Proof of Proposition \ref{prop1}.}}}
The proof uses arguments of the proof of \cite[Theorem 1]{CV}. 
 We need to show that $v=v^{\ell,\,\ell\lambar}_{m,\,k\theta}$ satisfies the defining relations of 
 $V(\bm{\xi}^{\ell,\,\ell\lambar}_{m,\,k\theta})$. The relations \eqref{cvs:dfr1} with $\mu=\ell\lambda+k\theta$ are clear from \eqref{mr:dfr1}. 
 We now prove the relations \eqref{cvs:dfr3}. Let $\alpha\in R^+$ and let $s_\alpha$
 denote the number of non-zero parts of $\xi(\alpha)={\xi}^{\ell,\,\ell\lambar}_{m,\,k\theta}(\alpha)$. 
 In the cases when
 (i) $r\geq\xi(\alpha)_1$, (ii) $r\leq \xi(\alpha)_{s_\alpha}$, and (iii) $\xi(\alpha)_{s_\alpha-1}>r>\xi(\alpha)_{s_\alpha}$, 
 the proof follows as in the proof of  \cite[Theorem 1]{CV} by using the relations \eqref{mr:dfr3}, \eqref{redrel} with $i=0$, and \eqref{mr:dfr5}-\eqref{mr:dfr7}.
 
We now prove the relations \eqref{cvs:dfr3} in the remaining case when $\xi(\alpha)_{1}> r\geq\xi(\alpha)_{s_\alpha-1}$.  
We observe from \eqref{xilmk} that this case is possible
only when $\alpha=\theta$ and $ m\leq 2k$. In this case, we have 
\beqs%\label{p:e1}
s_\theta=(\lambda|\theta)+2 \qquad \textup{and} \qquad \xi(\theta)_j=\begin{cases} \ell, & 1\leq j\leq s_\theta-2, \\ m, & j=s_\theta-1, \\2k-m, &j=s_\theta.
\end{cases}\eeqs
We observe that if  $(b_p)_{p\geq0}\in\mathbf{S}(r,s)$ is such that $b_p>0$ for some $p\geq s_\theta$, then by the relation \eqref{mr:dfr3} with $\alpha=\theta$,
 we have
$$\big((x^{-}_{\theta} \otimes 1)^{(b_0)}\cdots(x^{-}_{\theta}\otimes t^p)^{(b_p)}\cdots(x^{-}_{\theta} \otimes t^s)^{(b_s)}\big)\,v=0.$$
Hence we get
\beq\label{prp:eq1}
\big(\mathbf{x}^{-}_{\theta}(r,s)-\mathbf{x}^{-}_{\theta}(r,s)_{s_\theta}\big)\,v=0.
\eeq
If $(b_p)_{p\geq0}\in\mathbf{S}(r,s)_{s_\theta}$, then  $s=(s_\theta-1)b_{s_\theta-1}+\cdots+b_1$. 
We consider two cases. First suppose that $b_{s_\theta-1}>2k-m$. Then, using the relation \eqref{mr:dfr5} we have
\beq\label{prp:eq11}
\big((x^{-}_{\theta} \otimes 1)^{(b_0)}\cdots(x^{-}_{\theta} \otimes t^{s_\theta-2})^{(b_{s_\theta-2})}(x^{-}_{\theta} \otimes t^{s_\theta-1})^{(b_{s_\theta-1})}\big)\,v=0.
\eeq
%$\mathbf{x}^{-}_{\theta}(r,s)\,v=0$. 
In the second case, 
we have $b_{s_\theta-1}\leq 2k-m$.
If we prove that 
\beq\label{p:e2} b_{s_\theta-2}\geq2(k-b_{s_\theta-1})+1,\eeq
then \eqref{prp:eq11} would follow in this case from the relations \eqref{redrel} with $\alpha=\theta$ and  $i=b_{s_\theta-1}$.
This would give us $\mathbf{x}^{-}_{\theta}(r,s)_{s_\theta}\,v=0$. Using \eqref{prp:eq1},
we will have $\mathbf{x}^{-}_{\theta}(r,s)\,v=0$ which would complete the proof.

To prove \eqref{p:e2}, we observe that 
$$(s_\theta-1)b_{s_\theta-1}+(s_\theta-2)b_{s_\theta-2}+(s_\theta-3)(r-b_{s_\theta-1}-b_{s_\theta-2})\geq s\geq 1+r(\kay-1)+\sum_{j\geq\kay+1}\xi(\theta)_j,$$
which implies
$$2b_{s_\theta-1}+b_{s_\theta-2}\geq1+r(\kay-s_\theta+2)+\sum_{j\geq\kay+1}\xi(\theta)_j.$$
Since $r\geq \xi(\theta)_{s_\theta-1}$, we see that \eqref{p:e2} is immediate if $\kay\geq s_\theta-1$. If $\kay<s_\theta-1$, then 
$$2b_{s_\theta-1}+b_{s_\theta-2}\geq1+\sum_{s_\theta-1>j\geq\kay+1}(\xi(\theta)_j-r)+\xi(\theta)_{s_\theta-1}+\xi(\theta)_{s_\theta}\geq 2k+1,$$
where the last inequality is because $\xi(\theta)_j=\xi(\theta)_1 >r$ for all $1\leq j\leq {s_\theta-2}$.
%This completes the proof of the proposition.
\qed
%\begin{proposition}
%With hypothesis and notation as in Theorem \ref{MTCV}, the assignment $v^{\ell,\,\ell\lambar}_{m,\,k\theta}\mapsto v_{\bm{\xi}^{\ell,\,\ell\lambar}_{m,\,k\theta}}$ gives a surjective
%morphism from 
% $\mathbf{V}^{\ell,\,\ell\lambar}_{m,\,k\theta}$ onto $V(\bm{\xi}^{\ell,\,\ell\lambar}_{m,\,k\theta})$.
%\end{proposition}
%\begin{proof}
\subsubsection{{{Proof of Proposition \ref{prop2}.}}}
The following lemma is crucial in proving the proposition.
\begin{lemma}\label{lemCV}
Let $k, \ell, m\in\mathbb{N}$ be such that $\ell\geq m\geq k$ and $\lambda\in L^+$. 
%Let  $\lambar=(\lambda_1,\ldots,\lambda_p)\in ({L^+})^{p}$,  for $p\in\mathbb{N}$, and denote that $\lambda=\lambda_1+\cdots+\lambda_p$. 
Then, for all $\alpha\in R^+\setminus\{\theta\}$ and $1\leq i \leq k$, the relation
\beq\label{redrelforalpha}
\big(x^{-}_{\alpha}\otimes t^{(\lambda|\alpha)}\big)^{\langle (k-i)\theta,\, \alpha^{\vee} \rangle +1} \big(x^{-}_{\theta}\otimes t^{(\lambda|\theta)+1}\big)^i\,v^{\ell,\,\ell\lambar}_{m,\,k\theta}=0, 
\eeq is redundant
%in $\mathbf{V}^{\ell,\,\ell\lambar}_{m,\,k\theta}$, are redundant.
in the definition of $\mathbf{V}^{\ell,\,\ell\lambar}_{m,\,k\theta}$.
\end{lemma}
\begin{proof}
To prove the lemma, we show by considering two cases that the relation \eqref{redrelforalpha} follows from the relations \eqref{mr:dfr3} and \eqref{redrel} with $i=0$.\\

\noindent
{\bf Case (1).} {\em Suppose that $\alpha$ is either a long root or a short simple root.}\\

Since 
 $\alpha\in R^+\setminus\{\theta\}$, it is easy to see that $(\theta|\alpha)$ is either $0 \,\,\textup{or}\,\,1$.
If $(\theta|\alpha)=0$, then \eqref{redrelforalpha} follows from \eqref{mr:dfr3}. Assume that $(\theta|\alpha)=1$. We may also assume that $d_\alpha\in\{1,2\}$.
Indeed, if $d_\alpha=3$, then $\lieg$ is of type $G_2$ and  $\alpha$ is the short simple root of $\lieg$, hence $(\theta|\alpha)=0$.
We now have
\beq\label{thetaalpha}
(\theta|r_{\alpha}(\theta))=(\theta|\theta-d_\alpha\alpha)=2-d_\alpha\in\{0,1\}\qquad\textup{and}\qquad(\theta-\alpha), r_\alpha(\theta)\in R^+.\eeq 
Let $A_\alpha, B_\alpha\in\complex\setminus\{0\}$ be such that 
$[x^{-}_{r_{\alpha}(\theta)}, x^{-}_{\alpha}]=A_\alpha x^{-}_{\theta-(d_\alpha-1)\alpha}$ and $[x^{-}_{\theta-\alpha}, x^{-}_{\alpha}]=B_\alpha x^{-}_{\theta}.$
We denote $$f(A_\alpha, B_\alpha)=\begin{cases} A_\alpha=B_\alpha, & d_\alpha=1,\\
                                                           \frac{A_\alpha B_\alpha}{2}, & d_\alpha=2.
                                                                      \end{cases}$$                                                                      
Set $X=x^{-}_{r_{\alpha}(\theta)}\otimes t^{(\lambda|r_{\alpha}(\theta))+1}, Y=x^{-}_{\alpha}\otimes t^{(\lambda|\alpha)}, Z=x^{-}_{\theta}\otimes t^{(\lambda|\theta)+1}, 
\,\,\textup{and} \,\,v=v^{\ell,\,\ell\lambar}_{m,\,k\theta}$. 
We observe that 
\beq\label{eqqqqq}
[X, Z]=[Y, Z]=[[X, Y], Z]=0, \quad d_\alpha f(A_\alpha, B_\alpha)Z=\begin{cases}   [X, Y], & d_\alpha=1,\\
                                         [[X, Y], Y], & d_\alpha=2,
                           \end{cases} 
                           \qquad\textup{and}\quad X\,v=0,
\eeq
by \eqref{mr:dfr3} and \eqref{thetaalpha}.

%Let $i\in\mathbb{N}$ with $1\leq i\leq k$. 
By acting both sides to the relation $Y^{kd_\alpha+1}\,v=0$, which holds in this case in $\mathbf{V}^{\ell,\,\ell\lambar}_{m,\,k\theta}$ by 
the relation \eqref{redrel} with $i=0$,  with $X^{i}$, we get
\beq\label{cvs:e1}
X^{i}\,Y^{kd_\alpha+1}\,v=0.
\eeq
We claim that $\textup{for}\,\, \iprime\in\mathbb{Z}_{\geq 0},$
\beq\label{cvs:e2}
%\,\,\textup{the following relations hold in}\,\,\mathbf{V}^{\ell,\,\ell\lambar}_{m,\,k\theta}:\\
X^{(\iprime)}\,Y^{(k^\prime d_\alpha+1)}\,v
=(f(A_\alpha, B_\alpha))^{\iprime}Y^{((k^\prime-\iprime)d_\alpha+1)}\,Z^{(\iprime)}\,v,
\qquad \forall\,\,k^\prime\geq \iprime,
\eeq
where $\Bbb{X}^{(p)}$ denotes the divided power ${\Bbb{X}^p}/{p!}$. 
We observe from \eqref{cvs:e1} that the proof of \eqref{redrelforalpha} in this case  follows once we establish the claim. 

We now prove the claim \eqref{cvs:e2} by induction on $\iprime$.
For $\iprime=0$, there is nothing to prove. Now suppose that $\iprime\geq1$. Let $k^\prime\geq \iprime$. By the induction hypothesis, we have
\begin{equation}\label{cvs:e3}
%\begin{split}
X^{(\iprime-1)}\,Y^{(k^\prime d_\alpha+1)}\,v\\
=(f(A_\alpha, B_\alpha))^{\iprime-1}Y^{((k^\prime-\iprime+1)d_\alpha+1)}\,Z^{(\iprime-1)}\,v.
%\qquad \forall\,\,k^\prime\geq (i-1).
%\end{split}
\end{equation}
%for every $k^\prime\geq (\iprime-1)$. 
By acting both sides to \eqref{cvs:e3} with $\frac{1}{\iprime}X$,  we get
\begin{equation}\label{cvs:e4}
X^{(\iprime)}\,Y^{(k^\prime d_\alpha+1)}\,v=\frac{(f(A_\alpha, B_\alpha))^{\iprime-1}}{\iprime}
X\,Y^{((k^\prime-\iprime+1)d_\alpha+1)}\,Z^{(\iprime-1)}\,v.
\end{equation}
%for every $k^\prime\geq (\iprime-1)$. 
Since 
$X\,Z^{(\iprime-1)}\,v=
Z^{(\iprime-1)}\,X\,v=0$ by \eqref{eqqqqq}, 
we may replace the right hand side of \eqref{cvs:e4} by
\begin{equation}\label{cvs:e5}
\begin{split}
&
\frac{(f(A_\alpha, B_\alpha))^{\iprime-1}}{\iprime}[X, Y^{((k^\prime-\iprime+1)d_\alpha+1)}]\,Z^{(\iprime-1)}\,v\\
&=\frac{(f(A_\alpha, B_\alpha))^{\iprime-1}}{((k^\prime-\iprime+1)d_\alpha+1)!\iprime}\sum_{p=0}^{(k^\prime-\iprime+1)d_\alpha}Y^{p}\,[X, Y]\,Y^{(k^\prime-\iprime+1)d_\alpha-p}
\,Z^{(\iprime-1)}\,v.
%&=Y^{(k^\prime-\iprime+1)}\,Z^{(\iprime)}\,v.
\end{split}
\end{equation}

Suppose that $d_\alpha=1$. Then, since $[X, Y]=f(A_\alpha, B_\alpha)Z$ and $[Z, Y]=0$ by \eqref{eqqqqq}, it is easily checked that \eqref{cvs:e5} proves the claim \eqref{cvs:e2}.
Suppose that $d_\alpha=2$. Then, since $[[X, Y], Z]=0$ and $[X, Y]\,v=0$ by \eqref{mr:dfr3}, we have
\begin{equation}\label{cvs:e55}
\begin{split}
[X, Y]\,Y^{(k^\prime-\iprime+1)d_\alpha-p}\,Z^{(\iprime-1)}\,v&=[[X, Y], Y^{(k^\prime-\iprime+1)d_\alpha-p}]\,Z^{(\iprime-1)}\,v\\
&=\sum_{q=0}^{(k^\prime-\iprime+1)d_\alpha-p-1}Y^{q}\,[[X, Y], Y]\,Y^{(k^\prime-\iprime+1)d_\alpha-p-q-1}\,Z^{(\iprime-1)}\,v.
\end{split}
\end{equation}
Since $[[X, Y], Y]=2f(A_\alpha, B_\alpha)Z$ and $[Z, Y]=0$ by \eqref{eqqqqq}, the right hand side of \eqref{cvs:e55} simplifies to
$$2i^\prime f(A_\alpha, B_\alpha)((k^\prime-\iprime+1)d_\alpha-p)Y^{(k^\prime-\iprime+1)d_\alpha-p-1}\,Z^{(\iprime)}\,v.$$
Substituting this into \eqref{cvs:e5}, we get that the right hand side of \eqref{cvs:e4} is equal to
\begin{equation}\label{cvs:e555}
\begin{split}
&
\frac{2(f(A_\alpha, B_\alpha))^{\iprime}}{((k^\prime-\iprime+1)d_\alpha+1)!}\sum_{p=0}^{(k^\prime-\iprime+1)d_\alpha} ((k^\prime-\iprime+1)d_\alpha-p)Y^{(k^\prime-\iprime+1)d_\alpha-1}\,Z^{(\iprime)}\,v
\\&=\frac{2(f(A_\alpha, B_\alpha))^{\iprime}}{((k^\prime-\iprime+1)d_\alpha+1)!}\frac{((k^\prime-\iprime+1)d_\alpha)((k^\prime-\iprime+1)d_\alpha+1)}{2}Y^{(k^\prime-\iprime+1)d_\alpha-1}\,Z^{(\iprime)}\,v.
\end{split}
\end{equation}
This proves the claim \eqref{cvs:e2}.\\

\noindent
{\bf Case (2).} {\em Suppose that $\alpha$ is a short root.
}\\

We proceed by induction on $\textup{ht}\,\alpha$. If $\textup{ht}\,\alpha=1$, then \eqref{redrelforalpha} follows from case (1).
Assume that $\textup{ht}\,\alpha>1$. Then, since $\alpha$ is short, there exist a short root $\beta\in R^+$ and a root $\gamma\in R^+$ such that $\alpha=\beta+\gamma$.
Set $v_i=\big(x^{-}_{\theta}\otimes t^{(\lambda|\theta)+1}\big)^i\,v^{\ell,\,\ell\lambar}_{m,\,k\theta}$. 
By induction hypothesis, we have
\beq\label{eeeeeeeeeqqqqq1}
(x^{-}_{\beta}\otimes t^{(\lambda|\beta)}\big)^{\langle(k-i)\theta,\, \beta^\vee\rangle+1}\,v_i=0.
\eeq
If $\gamma$ is long (resp. short), then from case (1) (resp. induction hypothesis), we have
\beq\label{eeeeeeeeeqqqqq2}
(x^{-}_{\gamma}\otimes t^{(\lambda|\gamma)}\big)^{\langle(k-i)\theta,\, \gamma^\vee\rangle+1}\,v_i=0.
\eeq

It is easily checked that  the Lie subalgebra of $\lieg[t]$ generated by 
$\{x^{-}_{\beta}\otimes t^{(\lambda|\beta)}, x^{-}_{\gamma}\otimes t^{(\lambda|\gamma)}\}$ 
is isomorphic to the nilradical of the Borel subalgebra of the Lie algebra of type $A_2$ (resp. $B_2$, $G_2$), 
if $\gamma$ is short (resp. $\gamma$ is long and $d_\beta=2$, $\gamma$ is long and $
d_\beta=3$). Since
$$
\alpha^\vee=\begin{cases}\beta^\vee+d_\alpha\gamma^\vee,& \gamma \,\,\textup{is long},\\
\beta^\vee+\gamma^\vee,& \gamma \,\,\textup{is short},
                                         \end{cases}$$
we have from \eqref{eeeeeeeeeqqqqq1} and \eqref{eeeeeeeeeqqqqq2}
 (see \cite[Lemma 4.5]{Naoi1}) that
$$(x^{-}_{\alpha}\otimes t^{(\lambda|\alpha)}\big)^{\langle(k-i)\theta,\, \beta^\vee\rangle+d_\alpha\langle(k-i)\theta,\, \gamma^\vee\rangle+1}\,v_i=0, \qquad\textup{if}\,\,\gamma \,\,\textup{is long},$$
and $$(x^{-}_{\alpha}\otimes t^{(\lambda|\alpha)}\big)^{\langle(k-i)\theta,\, \beta^\vee\rangle+\langle(k-i)\theta,\, \gamma^\vee\rangle+1}\,v_i=0, \qquad\textup{if}\,\,\gamma \,\,\textup{is short},$$
which are equivalent to \eqref{redrelforalpha}.
Hence the lemma.
\end{proof}
We now prove Proposition \ref{prop2}.
Let $\bm{\xi}=\bm{\xi}^{\ell,\,\ell\lambar}_{m,\,k\theta}$. 
 We need to show that $v_{\bm{\xi}}$ satisfies the defining relations of  $\mathbf{V}^{\ell,\,\ell\lambar}_{m,\,k\theta}$.
 The relations  \eqref{mr:dfr1} are clear from \eqref{cvs:dfr1}. By taking $r=1, s=(\lambda+\theta|\alpha)$, and $\kay=(\lambda+\theta|\alpha)$ in \eqref{cvs:dfr3}, we obtain 
 \beq\label{pr:mr1}
 \big(x^{-}_{\alpha}\otimes t^{(\lambda+\theta|\alpha)}\big)\,v_{\bm{\xi}}=0, \quad \forall \,\, \alpha\in R^+,\eeq
 which proves the relations \eqref{mr:dfr3}.

To prove the relations \eqref{redrel}, using Lemma \ref{lemCV}, it is enough to prove that the following relations hold in $V(\bm{\xi})$
 \begin{gather}
\big(x^{-}_{\alpha}\otimes t^{(\lambda|\alpha)}\big)^{\langle k\theta,\, \alpha^{\vee} \rangle +1}\,v_{\bm{\xi}}=0, \quad \forall\,\,\alpha\in R^+,\label{enough:e1}\\
\big(x^{-}_{\theta}\otimes t^{(\lambda|\theta)}\big)^{2(k-i)+1} \big(x^{-}_{\theta}\otimes t^{(\lambda|\theta)+1}\big)^i\,v_{\bm{\xi}}=0, \quad\forall\,\,0\leq i\leq k\label{enough:e2}.
 \end{gather}
%Using \eqref{} and \eqref{pr:mr1}, the relations  \eqref{enough:e1} follow  from \eqref{cvs:dfr3}  by taking $r=\langle k\theta,\, \alpha^{\vee} \rangle +1, s=r(\lambda|\alpha)$, and $\kay=(\lambda|\alpha)$.
Using \eqref{pr:mr1}, the relations \eqref{enough:e1}  follow  from \eqref{cvs:dfr4}, by taking $r= \langle k\theta,\, \alpha^{\vee} \rangle +1$, 
$s=r(\lambda|\alpha)$, and $\kay=(\lambda|\alpha)$.
Using the relation $\big(x^{-}_{\theta}\otimes t^{(\lambda|\theta)+2}\big)\,v_{\bm{\xi}}=0$, 
the relations \eqref{enough:e2} (resp.  \eqref{mr:dfr5}, \eqref{mr:dfr6}) follow  from \eqref{cvs:dfr4} with $\alpha=\theta$,
by taking 
$r=2k-i+1$ (resp. $ r=2k-m+1$, $r=1$), $s=r(\lambda|\theta)+i$ (resp. $s=r(\lambda|\theta)+r$, $s=(\lambda|\theta)+1$), 
and $\kay=(\lambda|\theta)$  (resp. $\kay=(\lambda|\theta)+1$, $\kay=(\lambda|\theta)+1$). This completes the proof of the proposition.\qed
% \end{proof}
%\nocite{*}   %This command (when this line is uncommented)
%		will put all the entries in the database of references
%		into the body of the paper,  whether or not they are
%		cited in the text.    Without this only those cited 
%		will appear.

%%%%%%%%%%%%%%%%%%%%%%%%%%%%%%%%%%%%%%%%%%%%%%%%%%%%%%%%%%%%%%%%%%%%%%%%%%%%
%%%%%%%%%%%%%%% reformatting the citation number and number of the bibliography
%%%%%%%%	item can be achieved by the following:
%	\newcommand\citenumfont[1]{\textbf{#1}}
\bibliographystyle{bibsty-final-no-issn-isbn}%{plain}
\addcontentsline{toc}{section}{References}
%\ifthenelse{\equal{\finalized}{no}}{
%\bibliography{abbrev,references}}
%\input{tex_me.bbl}

\end{document}